\documentclass[11pt,a4paper]{article}
\usepackage{amsfonts,amsmath,amssymb,accents,amsthm}%,mathtools}
\usepackage{verbatim}
\usepackage{hyperref}
\usepackage[normalem]{ulem}%for \sout,\xout

\usepackage[notref,notcite]{}%{showkeys}

\usepackage{tikz} %also makes \color work
\newcommand{\dis}{\displaystyle}

\newcommand{\noi}{\noindent}
\newcommand{\halmos}{\rule{1ex}{1.4ex}}
\newcommand{\QED}{\nopagebreak{\hspace*{\fill}$\halmos$\medskip}}

\newcommand{\med}{\medskip}
\newcommand{\quand}{\quad\mbox{and}\quad}

\newtheoremstyle{mythm}% name
  {}%      Space above
  {}%      Space below
  {\itshape}%         Body font
  {}%         Indent amount (empty = no indent, \parindent = para indent)
  {\bfseries}% Thm head font
  {}%        Punctuation after thm head
  {.5em}%     Space after thm head: " " = normal interword space;
  {#1 #2 \thmnote{(#3)}}%  Thm head spec (can be left empty, meaning `normal')

\theoremstyle{mythm}
\newtheorem{theorem}{Theorem}[section]
\newtheorem{proposition}[theorem]{Proposition}
\newtheorem{lemma}[theorem]{Lemma}
\newtheorem{exercise}[theorem]{Exercise}
\newtheorem{corollary}[theorem]{Corollary}
\newtheorem{conjecture}[theorem]{Conjecture}

\newtheorem{counterex}[theorem]{Counterexample}
\newtheorem{remark}[theorem]{Remark}
\newtheorem{defi}[theorem]{Definition}
\newcommand{\bt}{\begin{theorem}}
\newcommand{\et}{\end{theorem}}
\newcommand{\bl}{\begin{lemma}}
\newcommand{\el}{\end{lemma}}
\newcommand{\bp}{\begin{proposition}}
\newcommand{\ep}{\end{proposition}}
\newcommand{\bcor}{\begin{corollary}}
\newcommand{\ecor}{\end{corollary}}
\newcommand{\br}{\begin{remark}\rm}
\newcommand{\er}{\end{remark}}
\newcommand{\bcon}{\begin{conjecture}}
\newcommand{\econ}{\end{conjecture}}
\newcommand{\bex}{\begin{exercise}}
\newcommand{\eex}{\end{exercise}}
\newcommand{\bcou}{\begin{counterex}}
\newcommand{\ecou}{\end{counterex}}

\newenvironment{Proof}[1][]{\noi\textbf{Proof#1.~}}{\QED}
\newcommand{\bpro}{\begin{Proof}}
\newcommand{\epro}{\end{Proof}}

\newcommand{\be}{\begin{equation}}
\newcommand{\ee}{\end{equation}}
\newcommand{\ba}{\begin{array}}
\newcommand{\ea}{\end{array}}
\newcommand{\bc}{\be\begin{array}{r@{\,}c@{\,}l}}
\newcommand{\ec}{\end{array}\ee}

\newcommand{\Ga}{\Gamma}
\newcommand{\de}{\delta}

\newcommand{\eps}{\varepsilon}

\newcommand{\sig}{\sigma}

\newcommand{\tet}{\theta}

\newcommand{\Di}{{\cal D}}

\newcommand{\R}{{\mathbb R}}
\newcommand{\N}{{\mathbb N}}
\newcommand{\Z}{{\mathbb Z}}

\newcommand{\E}{{\mathbb E}}
\renewcommand{\P}{{\mathbb P}}

\newcommand{\up}{\uparrow}
\newcommand{\down}{\downarrow}

\newcommand{\asto}[1]{\underset{{#1}\to\infty}{\longrightarrow}}
\newcommand{\Asto}[1]{\underset{{#1}\to\infty}{\Longrightarrow}}

\newcommand{\ov}{\overline}

\newcommand{\ffrac}[2]{{\textstyle\frac{{#1}}{{#2}}}}

\newcommand{\di}{\mathrm{d}}

\newcommand{\ha}{\ffrac{1}{2}}

\makeatletter
\let\@fnsymbol\@arabic
\makeatother

\begin{document}

%numbering formulas within sections
\makeatletter\@addtoreset{equation}{section}
\makeatother\def\theequation{\thesection.\arabic{equation}}

%alternative layout for enumerate lists.
\renewcommand{\labelenumi}{{\rm (\roman{enumi})}}
\renewcommand{\theenumi}{\roman{enumi}}

\title{Equilibrium interfaces of biased voter models}
\author{Rongfeng Sun \footnote{Department of Mathematics,
		National University of Singapore,
		10 Lower Kent Ridge Road, 119076 Singapore.
		Email: matsr@nus.edu.sg}
	\and Jan~M.~Swart \footnote{The Czech Academy of Sciences,
                Institute of Information Theory and Automation,
		Pod vod\'arenskou v\v e\v z\' i 4,
		18208 Prague 8,
		Czech Republic.
		Email: swart@utia.cas.cz}
	\and Jinjiong Yu \footnote{NYU-ECNU Institute of
		Mathematical Sciences at NYU Shanghai,
		3663 Zhongshan Road North,
		Shanghai 200062, China.
		Email: jinjiongyu@nyu.edu}
}

\date{\today}

\maketitle

\begin{abstract}\noi
A one-dimensional interacting particle system is said to exhibit interface
tightness if  \linebreak starting in an initial condition describing the interface between
two constant configurations of different types, the process modulo translations is positive
recurrent. In a biological setting, this describes two populations that do not
mix, and it is believed to be a common phenomenon in one-dimensional particle
systems. Interface tightness has been proved for voter models
satisfying a finite second moment condition on the rates.
We extend this to biased voter models.
%Our proof is an adaptation of a method of Sturm and Swart and
%is based on something akin to a Lyapunov function. This contrasts with earlier
%work which relied on duality.
Furthermore, we show that the distribution of the equilibrium interface for the biased voter model converges to that of the voter model when the bias parameter tends to zero. A key ingredient is an identity for the expected number of boundaries in the equilibrium voter model interface, which is of independent interest.

	%
	%Consider a one-dimensional interacting particle system in which each state can
	%be in two states and for which the constant configurations are invariant under
	%the dynamics. Let such a system be started in an initial state that decribes an
	%interface between two constant configurations. If such a system modulo
	%translations is positive recurrent, then one says that it exhibits interface
	%tightness. In particular, this says that the system spends a positive fraction
	%of time in the simplest possible interface and the two populations do not mix
	%in a strong sense.
	%
	%Interface tightness is believed to be a common phenomenon in
	%one-dimensional particle systems. Cox and Durrett introduced the concept and
	%proved it for long-range voter models under a third moment condition on the
	%rates, which was relaxed to a second moment condition by Belhaouari,
	%Mountford, and Valle. These papers depended heavily on duality of the voter
	%model to coalescing random walks. Later, Sturm and Swart gave a shorter proof
	%using a pseudo Lyapunov function. We show that the method of the latter paper
	%can be adapted to prove interface tightness also for biased voter models,
	%which are much harder to treat with duality methods.
\end{abstract}
\vspace{.5cm}

\noi
{\it MSC 2010.} Primary: 82C22; Secondary: 82C24, 82C41, 60K35.\newline
%82C22 Interacting particle systems
%82C24 Interface problems; diffusion-limited aggregation
%82C41 Dynamics of random walks, random surfaces, lattice animals etc.
%60K35 Interacting random processes; statistical mechanics type models,
%82C26   	Dynamic and nonequilibrium phase transitions (general)
%82B43   	Percolation
%82B26   	Phase transitions (general)
{\it Keywords.} Biased voter model, interface tightness, branching and
coalescing random walks.

{\setlength{\parskip}{-2pt}\tableofcontents}
%\newpage

\section{Introduction}
\subsection{Interface tightness}

One-dimensional \emph{biased voter models}, also known as one-dimensional
Williams-Bjerknes models \cite{S77,WB72}, are Markov processes
$(X_t)_{t\geq0}$ taking values in the space $\{0,1\}^\Z$ of infinite sequences
$x=(x(i))_{i\in\Z}$ of zeros and ones. They have several interpretations, one
of which is to model the dynamics of two biological populations. We call
$X_t(i)$ the type of the individual at site $i\in\Z$ at time
$t\geq0$. Let $\eps\in[0,1)$, and let $a:\Z\to[0,\infty)$ be a function such that $\sum_ka(k)<\infty$
and the continuous-time random walk that jumps from $i$ to $j$ with rate
$a(j-i)$ is irreducible. The dynamics of
$(X_t)_{t\geq0}$ are such that for each $i,j\in\Z$, at the times $t$ of a
Poisson point process with rate $a(j-i)$, if the type at $i$ just prior to
$t$ satisfies $X_{t-}(i)=1$, then $j$ adopts the type $1$; if $X_{t-}(i)=0$,
then with probability $1-\eps$, site $j$ adopts the type $0$, and with the
remaining probability $\eps$, $X_t(j)$ remains unchanged.

Somewhat more formally, we can define $(X_t)_{t\geq0}$ by specifying its
generator. For $x\in\{0,1\}^\Z$ and $i_1,\ldots,i_n\in\Z$, write
$x(i_1,\ldots,i_n):=(x(i_1),\ldots,x(i_n))\in\{0,1\}^n$. We also use the convention
of writing elements of $\{0,1\}^n$ as words consisting of the letters $0$ and
$1$, i.e., instead of $(1,0)$ we simply write $10$ and similarly for longer
sequences. With this notation, the generator of the biased voter model we
have just described is given by
\bc\label{bvmgen}
\dis G^\eps f(x)
&=&\dis\sum_{i,j}a(j-i)1_{\{x(i,j)=10\}}\big\{f(x+e_j)-f(x)\big\}\\[5pt]
&&\dis+(1-\eps)\sum_{i,j}a(j-i)1_{\{x(i,j)=01\}}\big\{f(x-e_j)-f(x)\big\},
\ec
where $e_i(j)\!:=\!1_{\{i=j\}}$. We call $\eps$ the \emph{bias} and $a$ the
\emph{underlying random walk kernel}. In particular, for $\eps=0$, we obtain
a normal (unbiased) voter model, in which the types $0$ and $1$ play symmetric
roles. By contrast, for $\eps>0$, the $1$'s replace $0$'s at a faster rate
than the other way round, i.e., there is a bias in favor of the $1$'s. To indicate the bias
parameter, $(X_t)_{t\geq0}$ will be denoted as $(X^\eps_t)_{t\geq0}$
hereafter.

Define
\bc\label{eq:S01}
\dis S^{01}_{\rm int}&:=&\dis\big\{x\in\{0,1\}^\Z:
\lim_{i\to-\infty}x(i)=0,\ \lim_{i\to\infty}x(i)=1\big\},\\[5pt]
\dis S^{10}_{\rm int}&:=&\dis\big\{x\in\{0,1\}^\Z:
\lim_{i\to-\infty}x(i)=1,\ \lim_{i\to\infty}x(i)=0\big\},
\ec
which denote the sets of states describing the interface between an infinite cluster
of 1's and an infinite cluster of 0's. 
We can also define $S^{00}_{\rm int}$ and $S^{11}_{\rm int}$ analogously by changing the limiting behavior at $\pm\infty$. 
If $\sum_ka(k)|k|<\infty$ and
$0\leq\eps<1$, then it is not hard to see that $X^\eps_0\in S^{01}_{\rm int}$
implies that $X^\eps_t\in S^{01}_{\rm int}$ for all $t\geq 0$, a.s., and
similarly for $S^{10}_{\rm int}$. (For unbiased voter models, this is proved
in \cite{BMSV06}. The proof in the biased case is the same.) We will be
interested in studying the long-time behavior of the interface of
$(X^\eps_t)_{t\geq0}$.

We call two configurations $x,y\in\{0,1\}^\Z$ \emph{equivalent}, denoted by
$x\sim y$, if one is a translation of the other, i.e., there exists some
$k\in\Z$ such that $x(i)=y(i+k)$ $(i\in\Z)$. We let $\ov x$ denote the
equivalence class containing $x$ and write
\be
\ov S^{01}_{\rm int}:=\{\ov x:x\in S^{01}_{\rm int}\}
\quand
\ov S^{10}_{\rm int}:=\{\ov x:x\in S^{10}_{\rm int}\}.
\ee
Note that $S^{01}_{\rm int},\ov S^{01}_{\rm int}$, $S^{10}_{\rm int}$ and $\ov S^{10}_{\rm int}$ are countable sets.
Since our rates are translation invariant, the \emph{process modulo
  translations} $(\ov{X}^\eps_t)_{t\geq0}$ is itself a Markov process.
For non-nearest neighbor kernels $a$, this Markov process is
  irreducible; see Lemma~\ref{L:irre} below. From now on, we restrict
  ourselves to the non-nearest neighbor case.
%
%However, even if $a$ is irreducible, $(\ov{X}^\eps_t)_{t\geq0}$ being
%  irreducible is not entirely trivial, as a counterexample arises when $a$ is
%  nearest neighbor.  But as will be shown by Lemma~\ref{L:irre} in
%  Subsection~\ref{S:tightout}, excluding this case the translated process
%  $(\ov{X}^\eps_t)_{t\geq0}$ is always irreducible, and thus we will restrict
%  ourselves in non-nearest neighbor cases hereafter.
%
We adopt the following definition from \cite{CD95}.
\begin{defi}[Interface tightness]\label{D:int}
We say that $(X^\eps_t)_{t\geq0}$ exhibits interface tightness on
$S^{01}_{\rm int}$ $($resp.\ $S^{10}_{\rm int}$$)$ if
$(\ov X^{\eps}_t)_{t\geq 0}$ is positive recurrent on $\ov S^{01}_{\rm int}$
$($resp.\ $\ov S^{10}_{\rm int}$$)$.
\end{defi}
\noi
Note that because of the bias, if $a$ is not symmetric, then there is no
obvious symmetry telling us that interface tightness on $S^{01}_{\rm int}$
implies interface tightness on $S^{10}_{\rm int}$ or vice versa.

In the unbiased case $\eps=0$, Cox and Durrett \cite{CD95} proved that
interface tightness holds if $\sum_{k}a(k)|k|^3<\infty$. This was relaxed to
$\sum_{k}a(k)k^2<\infty$ by Belhaouari, Mountford and Valle \cite{BMV07}, who
moreover showed that the second moment condition is optimal. Our first main result is the following theorem that extends
this to biased voter models, where the optimal condition in the biased case
turns out to be even weaker.

\bt[Interface tightness for biased voter models]
%Assume\label{T:tight} that $\eps\in[0,1)$ and that the kernel $a$ is
%non-nearest neighbor \Jan{and irreducible. If $\sum_{k<0}a(k)k^2<\infty$ and
%$\sum_{k>0}a(k)k<\infty$, then the biased voter model
%$(X^\eps_t)_{t\geq0}$ exhibits interface tightness on $S^{01}_{\rm int}$.
%Similarly, $\sum_{k<0}a(k)|k|<\infty$ and
%$\sum_{k>0}a(k)k^2<\infty$ imply interface tightness on $S^{10}_{\rm int}$.}
Assume\label{T:tight} that $\eps\in[0,1)$ and that the kernel $a$ is
non-nearest neighbor, irreducible and satisfies $\sum_ka(k)k^2<\infty$. Then
the biased voter model $(X^\eps_t)_{t\geq0}$ exhibits interface tightness on
$S^{01}_{\rm int}$ and $S^{10}_{\rm int}$. If $\eps>0$ and the moment
condition on the kernel is relaxed to $\sum_{k<0}a(k)k^2<\infty$ and
$\sum_{k>0}a(k)k<\infty$, then interface tightness on $S^{01}_{\rm int}$ still
holds.
\et

\noi
To see heuristically why for interface tightness on $S^{01}_{\rm int}$,
a finite first moment condition in the positive direction suffices, we
observe that $a(k)$ with large positive $k$ govern $0$'s that are created deep into the
territory of the $1$'s. Such $0$'s do not survive long due to the bias. On the
other hand, $a(k)$ with large negative $k$ govern $1$'s that are created deep into the
territory of the $0$'s. These $1$'s have a positive probability of surviving and
then giving birth to $1$'s even further away. This explains heuristically why
one needs to impose a stronger moment condition in the negative direction.

The proof of interface tightness for the voter model in \cite{CD95,BMV07} relied heavily on the
well-known duality of the voter model to coalescing random walks. A biased voter
model also has a dual, which is a system of coalescing random walks that
moreover branch with rate $\eps$. Because of the branching, this dual process
is much harder to control than in the unbiased case. In view of this, we were
unable to apply the methods of \cite{CD95,BMV07} but instead adapted a method
of \cite{SS08b}, who provided a short proof of interface tightness for
unbiased voter models using generator calculations and a Lyapunov like function.
Our key observation is that this function admits a generalization to biased voter models, and the proof of
interface tightness can then be adapted accordingly. However, as we will see in the proof of
Theorem~\ref{T:tight}, the calculations are considerably more complicated in the biased case.

We remark that interface tightness is a fairly common phenomenon among
  many one-dimensional interacting particle systems. Other models for
    which this has been proved include one-dimensional two-type contact
  processes with strong bias (in the sense that one type can overtake the
    other, but not vice versa) \cite{AMPV10}, or no bias (no type can
    infect a site occupied by the other type) \cite{Val10, MV16}, as well as
  one-dimensional asymmetric exclusion processes that admit so-called blocking
  measures (see e.g.~\cite{BM02, BLM02}).

%Besides the voter model, %it has also been proved for one-dimensional multi-type contact processes, see e.g., \cite{AMPV10} and \cite{Val10}, where the proofs rely on duality, renormalization and percolation techniques. {\yu one-dimensional multi-type contact processes are also examples with a Glauber-type dynamics having this property, see e.g.,\ \cite{AMPV10} and \cite{Val10}, where the proofs rely on duality, renormalization and percolation techniques. Furthermore, based on the interface tightness result, Mountford and Valesin \cite{MV16} showed that the scaling limit of the interface is a Brownian motion. Examples with a Kawasaki-type dynamics include one-dimensional exclusion processes for which the underlying random walk kernel $a$ has nonzero mean, see e.g.\ \cite{BLM02} and \cite{BM02}.
%When the underlying random walk kernel $a$ has positive mean and finite range, Bramson and Mountford \cite{BM02} showed existence of the invariant measure of the exclusion process on a suitable subset of $S^{01}_{\rm int}$, and thus interface tightness.
%In \cite{BLM02}, this result is shown while $a$ being finite range is relaxed to $\sum_{k<0}a(k)k^2<\infty$, but two other conditions on $a$ are necessarily imposed.
%In contrast, if $\sum_{k<0}a(k)k^2=\infty$, then it is shown in the same paper that no such invariant measure exists.}

\subsection{Convergence of the equilibrium interface}\label{S:intrnucon}
%Continuity as a function of the bias}

Theorem \ref{T:tight} and Lemma~\ref{L:irre} imply that, modulo translations, the biased voter
model on $S^{01}_{\rm int}$ is
an irreducible and positive recurrent countable state Markov chain, and hence has a unique invariant law and is
ergodic. In particular, if $\sum_ka(k)k^2<\infty$, then for each $\eps\geq 0$, there is a unique
invariant law $\ov{\nu}_\eps$ on $\ov S^{01}_{\rm int}$. It is a natural question whether
$\ov{\nu}_\eps$ converges to the unique invariant law $\ov{\nu}_0$ for the voter model.
Our second main result answers this question affirmatively.

\bt[Continuity of the invariant law]\label{T:nucont}
Assume that the kernel $a(\cdot)$ is non-nearest neighbor, irreducible, and satisfies
$\sum_ka(k)k^2<\infty$. Then as $\eps\down 0$, the invariant law
$\ov{\nu}_\eps$ converges weakly to $\ov{\nu}_0$ with respect to the discrete
topology on $\ov{S}^{01}_{\rm int}$.
\et

\noi Recall that $\ov{S}^{01}_{\rm int}$ is the set of equivalence classes of
elements of $S^{01}_{\rm int}$ that are equal modulo translations. It will also be convenient to
choose a representative from each equivalence class by shifting the leftmost one to the origin. Since each equivalence class $\ov x\in\ov
S^{01}_{\rm int}$ contains a unique element $x$ in the set
\be\label{hatS'}
\hat{S}^{01}_{\rm int}:=\{x\in S^{01}_{\rm int}:x(i)=0
\mbox{ for all $i<0$ and }x(0)=1\},
\ee
we can identify $\ov S^{01}_{\rm int}$ with $\hat S^{01}_{\rm int}$. Under
this identification, $\ov{\nu}_\eps$ on $\ov{S}^{01}_{\rm int}$ uniquely
determines a probability measure $\nu_\eps$ on $\{0,1\}^\Z$ that is supported
on $\hat{S}^{01}_{\rm int}$.

The proof of Theorem \ref{T:nucont} turns out to be much more delicate than the result may suggest.
The difficulty lies in the choice of the discrete topology on $\ov{S}^{01}_{\rm int}$. In particular,
Theorem \ref{T:nucont} implies that the length of the equilibrium interface under $\ov{\nu}_\eps$
is tight as $\eps\downarrow 0$, and if started at the {\em heaviside state} $x_0$ with
\be\label{heavy}
x_0(i)=0\quad
{\rm for}~i<0\qquad{\rm and}\qquad
x_0(i)=1\quad{\rm for}~i\geq0,
\ee
then the time it takes to return to the state $x_0$ is also tight. Such uniform control in $\eps$ as $\eps\downarrow 0$
turns out to be difficult to obtain. We get around this difficulty by first proving the weak convergence of $\nu_\eps$
to $\nu_0$ under the product topology on $\{0, 1\}^\Z$. 
%\Swa{\sout{where $\nu_\eps$ and $\nu_0$ are induced by $\ov{\nu}_\eps$
%and $\ov{\nu_0}$ by choosing a representative in $\hat{S}^{01}_{\rm int}$ as defined above.}}
%$\ov{\nu}_\eps$ to $\ov{\nu}_0$ under the product topology, where we identify $\ov{S}^{01}_{\rm int}$ with a subset of the product space $\{0,1\}^\N$ by shifting the leftmost $1$ of the interface to the origin.
To strengthen this to convergence under the discrete topology, we
prove, for the unbiased equilibrium interface, an exact formula for
the expectation of a quantity involving the number of
  \emph{$k$-boundaries}, i.e., pairs of a $0$ and a $1$ at distance $k$ apart,
  as well as a matching one-sided bound in the biased case (see
  Proposition~\ref{P:equi} and Lemma~\ref{L:Ibdd}). These results
are also of independent interest.

\subsection{Relation to the Brownian net}\label{S:BN}

Theorem \ref{T:nucont} is in fact motivated by studies of branching-coalescing
random walks and their convergence to the Brownian net under weak branching. Let
us now explain this connection.

Similar to the well-known duality between the voter model and coalescing random walks, %\Swa{\sout{the}}
biased voter models are dual to systems of branching-coalescing random walks,
with the bias $\eps$ being the branching
rate of the random walks. While in \cite{CD95,BMV07}, coalescing random walks were used to
prove interface tightness, our motivation is the other way
around: we aim to use interface tightness as a tool to study the dual branching-coalescing random walks. More precisely, the present
paper arises out of the problem of proving convergence of rescaled
branching-coalescing random walks with weak branching to a continuum object called the
\emph{Brownian net}~\cite{SS08c}.

Let us first recall the graphical representation of (biased) voter
models. Plot space horizontally and time vertically, and for each $i,j\in\Z$, at the times $t$ of an independent Poisson point
process with intensity $(1-\eps)a(j-i)$, draw an arrow from $(i,t)$ to
$(j,t)$. We call such arrows \emph{resampling arrows}. Also, for each
$i,j\in\Z$, at the times $t$ of an independent Poisson process with intensity
$\eps a(j-i)$, draw a different type of arrow (e.g., with a different color)
from $(i,t)$ to $(j,t)$. We call such arrows \emph{selection arrows}.

It is well-known that a voter model can be constructed in such a way
that starting from the initial state, at each time $t$ where there is
a resampling arrow from $(i,t)$ to $(j,t)$, the site $j$ adopts the type of
site $i$. To get a biased voter model one also adds the selection arrows,
which are similar to resampling arrows, except that they only have an effect
when the site $i$ is of type 1.

To see the duality, we construct a system of coalescing random walks evolving
backwards in time as follows: let the backward random walk at site $j$ jump to
$i$ when it meets a resampling arrow from $i$ to $j$, and let it coalesce
with the random walk at site $i$ if there is one.
A system of backward branching-coalescing random walks can be obtained
by moreover allowing the coalescing random walks to branch at selection arrows.
That is, when a random walk at $j$ meets a selection arrow from $i$ to $j$, let
it branch into two walks located at $i$ and $j$, respectively.
The duality goes as follows. Let $A$ and $B$ be two sets of integers.
For the biased voter model starting from the state being $1$ only on $A$ at time $0$,
the set of $1$'s at time $t$ has nonempty intersection with $B$ if and only if for the
backward branching-coalescing random walks starting from $B$ at time $t$, there is at least
one walk in $A$ at time $0$.

It is shown in \cite{FINR04} that for coalescing nearest neighbor random walks
in the space-time plane, the diffusively rescaled system converges to the
so-called \emph{Brownian web}, which, loosely speaking, is a collection of
coalescing Brownian motions starting from every space-time point (see also \cite{SSS17}
for a survey on the Brownian web, Brownian net and related topics). Later, in
\cite{NRS05}, this result was extended to general coalescing random walks with
a finite fifth moment. This condition was then relaxed to a finite $(3+\eta)$-th
moment by Belhaouari et al.\ \cite{BMSV06}. It was observed in the same article
(see \cite[Theorem~1.2]{BMSV06}) that to verify tightness for rescaled systems
of random walks, it suffices to show that for the dual voter model starting
from the heaviside state $x_0$ (cf.~\ref{heavy}), the trajectories of the left and the right
interface boundaries converge to the same Brownian
motion. In the biased case, Sun and Swart \cite{SS08c} showed that systems
of branching-coalescing nearest neighbor random walks converge to the
Brownian net as the branching rate $\eps$ tends to zero and space and time are
diffusively rescaled with the same $\eps$. Just as in the unbiased case,
in order to extend the result of \cite{SS08c} to non-nearest neighbor random
walks, tightness can be established by showing that the left and right boundaries of the
dual biased voter model interface converge to the same Brownian motion, which
in contrast to the unbiased case, now has a drift. This, however, remains a challenge. Our Theorem
\ref{T:nucont} on the convergence in law of the equilibrium biased voter model interface, as the bias parameter $\eps\downarrow 0$,
can be regarded as a first step in this direction.

%Nevertheless, a first step in this direction is to prove that interface tightness holds in a way that is uniform as the bias parameter $\eps$ decreases to zero, which led to our Theorem \ref{T:nucont}.
%{\yu Following this general strategy, we can further show in \cite{?} weak convergence of a measure-valued process associated to the biased voter model, where the interface converges in a sense to the desired drifted Brownian motion.
%In essence, the measure-valued process discounts those isolated rare $0$'s and $1$'s in the diffusive scaling limit, which plays a key role in the proof.
%On the other hand, weak convergence of the left and right boundaries of the interface, however, remains a challenge.
%}
%{has been moved to the end of Section 1.1} We also remark that based on the interface tightness result obtained in \cite{Val10}, Mountford and Valesin \cite{MV16} showed that the interface of the multi-type contact process has Brownian motion as its scaling limit.

\subsection{Some open problems}
Before closing the introduction, we list some open problems in two directions.
First, what if we allow different kernels for resampling and selection?
%while believed to be a common phenomenon, interface tightness has so far been proved only
%for a few systems, including (biased) voter models and multi-type contact processes.
%Staying within the class of biased voter models, since there are
%two types of arrows, it would be meaningful to consider models with different kernels for resampling and selection arrows.
It seems plausible that for such models, Theorems~\ref{T:tight} and \ref{T:nucont} should still hold if both
kernels satisfy our moment assumptions. But unfortunately, our methods break
down if the kernels are different. Another class of models for which the question of  interface
tightness remains open are voter models with weak heterozygosity selection, such as
the rebellious voter model introduced in \cite{SS08a}.

Another direction of further research concerns the diffusive scaling limits of the biased
voter model as the bias parameter $\eps\downarrow 0$, as well as its connection with the Brownian net.
%is concerned with questions related to the models studied in Theorems~\ref{T:tight} and \ref{T:nucont}.
In particular, as already pointed out in the last section,
sending $\eps$ to zero and at the same time rescaling space by $\eps$ and time by
$\eps^2$, one would like to show that the interface converges to a drifted Brownian
motion, and then use this to prove convergence of the dual system of
branching-coalescing random walks to the Brownian~net.
A weaker formulation would be to consider the measure-valued process associated with the biased voter model (i.e., the counting measure of $1$'s in the baised voter configuration), and show that under diffusive scaling, it converges to a measure-valued process with a sharp interface between the regions with density either $0$ or $1$ with respect to the Lebesgue measure, such that the interface evolves as a drifted Brownian motion. For the voter model, such a result was established in \cite{AS11}.
% with boundary being the drifted Brownian motion. We remark that in the limit the counting measure would naturally discount isolated $0$'s and $1$'s of the biased voter model interface, and therefore this convergence is a weaker version compared with the convergence of the boundaries.

The rest of the paper is devoted to proofs. We prove Theorem~\ref{T:tight}
in Section~\ref{S:tight} and Theorem~\ref{T:nucont} in Section~\ref{S:nucont}.

\section{Interface tightness} \label{S:tight}

\subsection{Elementary observations and outline}\label{S:tightout}

Before moving to the proof of Theorem~\ref{T:tight}, we present the following elementary lemma that shows the irreducibility of non-nearest neighbor biased voter models.
\bl[Irreducibility of the biased voter model]
Assume\label{L:irre} that the kernel $a$ is non-nearest neighbor {\rm(}i.e.,
$a(k)>0$ for some $|k|\geq 2${\rm)}, irreducible and satisfies $\sum_ka(k)|k|<\infty$.
Then the biased voter model $(X^\eps_t)_{t\geq0}$ on $S^{01}_{\rm int}$ is irreducible.
\el
\bpro
As proved in \cite{BMSV06} and also in Lemma~\ref{L:EIeps10} below, the
condition $\sum_{k}a(k)|k|<\infty$ guarantees that the continuous-time Markov
chain on $S^{01}_{\rm int}$ is well-defined and nonexplosive.

We will show that for any configuration $x\in S^{01}_{\rm int}$, 1) for
  each $i\in\Z$, there is a positive probability to reach the heaviside state
  $x_i:=1_{\{i,i+1,\ldots\}}$ from $x$, and 2) there exists some $i\in\Z$ such
  that there is a positive probability to reach $x$ from $x_i$. The first
statement actually holds for any irreducible $a$. For such $a$, there exist
$k_{\rm r},k_{\rm l}>0$ such that $a(k_{\rm r}),a(-k_{\rm l})>0$. Since
$a(k_{\rm r})>0$ (resp.~$a(-k_{\rm l})>0)$, with positive probability the
leftmost $1$ (resp.~the rightmost $0$) can change into a $0$ (resp.~$1$)
while all other sites remain unchanged. In this way, $x_i$ can be
  reached for any $i\in\Z$.
%each $i\in\Z$ and $x\in S^{01}_{\rm int}$, $x_i$ can be reached from $x$.}

For the second statement, let $k$ be such that $a(k)>0$ for some $|k|\geq
  2$. We will prove that if $k<0$ (resp.\ $k>0$), then the interface can
  always be expanded by $1$ unit at the right (resp.\ left) boundary without
  changing the values at other sites, which implies statement 2). By symmetry,
  it suffices to consider the case $k<0$.
Suppose that the right boundary of $X^\eps_t$ is at site $i$, namely $X^\eps_t(i)=0$ and $X^\eps_t(j)=1$ for all $j>i$.
Then we will construct infections to show that at a later time $s$, with
  positive probability, it may happen that $X^\eps_s(i,i+1)=00$ (similarly, it may happen that $X^\eps_s(i,i+1)=10$), while the values of other sites of $X^\eps_t$ and $X^\eps_s$ are the same.
Indeed, since $a$ is irreducible, a path $\pi$ of infections from site $i$ to $i+1$ can be
constructed where moreover the path first does right jumps, and then left jumps.
Recall that the value of site $i$ is $0$.
Thus, the value of site $i+1$ is altered to $0$ and all sites on the left of $i$ remains unchanged.
Now by applying left infections of size $k$, one can consecutively alter the values back to $1$ for those sites on the right of $i+1$ that were infected by $\pi$.
After such infections, we have $X^\eps_s(i,i+1)=00$ while other sites remain unchanged.
To see that the value of site $i$ can also be altered to $1$, simply note that
site $i+|k|$ can infect site $i$ by a left infection of size $k$.
Furthermore, we have $X^\eps_s(i+|k|)=1$ since $|k|\geq 2$.
This completes the proof.
\epro

The irreducibility of the biased voter model $(X^\eps_t)_{t\geq0}$ immediately
implies the irreducibility of the process modulo translations $(\ov{X}^\eps_t)_{t\geq0}$.
Also note that $S^{01}_{\rm int}$, and hence $\ov{S}^{01}_{\rm int}$, is countable.
Therefore, by Definition~\ref{D:int}, $(X^\eps_t)_{t\geq0}$ exhibits interface tightness if and only if there exists an invariant probability measure for $(\ov{X}^\eps_t)_{t\geq0}$. %{\yu \sout{, where the latter implies that the interface being large at a fixed time has a small probability.
%More precisely,}
In particular, if $L:S^{01}_{\rm int}\to\N$, defined by
\be\label{intlen}
L(x):=\max\{i:x(i)=0\}-\min\{i:x(i)=1\}+1,
\ee
denotes the interface length, then interface tightness is equivalent to the family of random variables
$\big(L(X^\eps_t)\big)_{t\geq0}$ being tight. By definition, an
  \emph{inversion} is a pair $(i,j)$ such that $j<i$ and $x(j,i)=10$.
We let $h:S^{01}_{\rm int}\to\R$ denote the function counting the number of
inversions
\be\label{numinv}
h(x):=\big|\{(j,i):j<i,\ x(j,i)=10\}\big|,
\ee
where $|\cdot|$ denotes the cardinality of a set.
It is easy to see that $h(x)\geq L(x)-1$ for any configuration $x$.
In the proofs of \cite{CD95,BMV07}, the number of inversions $h(x)$ plays a key role, where duality is used to show that this quantity cannot be too big.
In \cite{SS08b}, $h$ is also a key ingredient, playing a role similar to a
Lyapunov function as in Foster's theorem (see e.g.,~\cite[Theorem 2.6.4]{MPW17}).
More precisely, it is shown there
that if interface tightness does not hold, then over sufficiently long time
intervals, $h$ would have to decrease on average more than it increases,
contradicting the fact that $h\geq 0$.

In the rest of this section, we will adapt the method of \cite{SS08b} to show
Theorem~\ref{T:tight}. In the biased case, we need a different function
  from $h$, which we call the {\em weighted number of inversions}.
In Subsection~\ref{S:Thm1} , we will state three necessary lemmas and then prove Theorem~\ref{T:tight}. We then prove those three lemmas in
Subsection~\ref{S:Gh}, Subsection~\ref{S:nonexp} and Subsection~\ref{S:intgr}, respectively.

\subsection{Proof of Theorem \ref{T:tight}}\label{S:Thm1}

Our key observation here is that
for the biased voter model, the number of inversions $h$ in \eqref{numinv} can be replaced
by $h^\eps:S^{01}_{\rm int}\to\R$ with
%the number of inversions $h$ in \eqref{numinv} can be {\yu\sout{generalized} adapted} to the biased case as follows. Let
%$h^\eps:S^{01}_{\rm int}\to\R$ be defined by
\be\label{hdef2}
h^\eps(x):=\sum_{i>j}(1-\eps)^{\sum_{n<j}x(n)}1_{\{x(j,i)=10\}}.
\ee
By numbering the $1$'s in $x$ from left to right, we see that $h^\eps(x)$ is
a weighted number of inversions, where inversions
involving the $j$-th $1$ carry weight $(1-\eps)^{j-1}$.
It is clear that for every configuration $x\in S^{01}_{\rm int}$, $h^0(x)$ agrees with the number of inversions $h(x)$ given in (\ref{numinv}), and %{\yu\sout{$h^{\eps}(x)\to h^{0}(x)$} }
$h^{\eps}(x)\uparrow h^{0}(x)$ as $\eps\downarrow0$.

To prove Theorem~\ref{T:tight}, we need three lemmas that generalize Lemma~2, Lemma~3 and Proposition~4 of \cite{SS08b} to the biased voter model.
To state the first lemma that gives an expression for the action of the generator from  \eqref{bvmgen} on $h^\eps$, we introduce the following notation.

By definition, a \emph{$k$-boundary} is a pair $(i,i+k)$ such that
$x(i)\neq x(i+k)$. For $k\in\Z$, let $I_k:S^{01}_{\rm int}\to\N$ be
the function counting the number of $k$-boundaries
\be\label{Idef}
I_k(x):=\big|\{i:x(i)\neq x(i+k)\}\big|.
\ee
\bl[Generator calculations]\label{L:Geps}
Under the assumption of Theorem~\ref{T:tight}, we have that for any $\eps\in[0,1)$ and $x\in S^{01}_{\rm int}$,
\be\label{Gh}
G^\eps h^\eps(x) = \sum_{k} a(k) \big(\ha k^2-\eps R^\eps_k(x)\big) - \ha\sum_k a(k) I_{k}(x),
\ee
%[Writing $G^\eps h^\eps$ in this way, if $\sum_{k}a(k)k^2=\infty$, then we would have a meaningless form $G^\eps h^\eps=\infty-\infty$. Shall we remark on this?]
where the generator $G^\eps$ is defined in \eqref{bvmgen}, and the term
$R^\eps_k(x)\geq0$ is given by $R^\eps_k(x):=0$ for $k=-1,0,1$ and
\be\label{Rdef}
R^\eps_k(x):=\left\{\ba{ll}
\dis\sum_{i}\sum_{n=1}^{k-1}(1-\eps)^{\sum_{j<i}x(j)}(k-n)1_{\{x(i-n,i)=01\}}&\dis(k>1),\\ [5pt]
\dis\sum_{i}\sum_{n=1}^{|k|-1}(1-\eps)^{\sum_{j<i}x(j)}(|k|-n)1_{\{x(i,i+n)=10\}}&\dis (k<-1).
\ea\right.
\ee
%Moreover, for $\eps\in[0,1)$ we have {\rm [to be deleted.]}
%\be\label{Ghupp}
%G^\eps h^\eps(x) \leq \ha\sum_{k}a(k)k^2-\ha\sum_{k}a(k)I_k(x).
%\ee
Moreover, for $\eps\in(0,1)$, we have
\be \label{Ghupp2}
G^\eps h^\eps(x) \leq \ha\sum_{k<0}a(k)k^2+\eps^{-1}\sum_{k>0}a(k)k
-\ha\sum_{k}a(k)I_k(x).
\ee
\el
\br
In the unbiased case $\eps=0$, \eqref{Gh} reduces to
\be\label{G0}
G^0 h^0(x) = \ha\sum_{k} a(k)\big(k^2-I_k(x)\big),
\ee
which agrees with \cite[Lemma~2]{SS08b}.
Since $R^\eps_k\geq0$ for all $k\neq0$, \eqref{Gh}  shows that $G^0h^0$ is an upper bound of $G^\eps h^\eps$ uniformly in $\eps$ in the sense that $G^\eps h^\eps(x)\leq G^0h^0(x)$.
\er

\br
In a sense, the best motivation we can give on why the weighted number of
inversions $h^\eps$ as in (\ref{hdef2}) is the ``right'' function to look at
is formula (\ref{Gh}). In the nearest-neighbor case
$a(-1)=a(1)=\ha$ and $a(k)=0$ for all $k\neq -1,1$, formula (\ref{Gh}) reduces
to $G^\eps h^\eps(x)=\ha(1-I_1(x))$. In particular, $G^\eps h^\eps(x)=-1$ if
$I_1(x)=3$, which can be used to prove (compare \cite[Lemma~12]{SS15}) that in
the nearest-neighbor case, starting from an initial state with $I_1(x)=3$,
$h^\eps(x)$ is the expected time before the system reaches a heaviside
state. These observations originally motivated us to define $h^\eps$ as in
(\ref{hdef2}).
\er

\noi We need two more lemmas.
Recall that $(X^\eps_t)_{t\geq0}$ denotes the biased voter model with bias $\eps$.
\bl [Nonnegative expectation]\label{L:nonexp}
Let the biased voter model $(X^\eps_t)_{t\geq0}$ start from a fixed configuration $X^\eps_0=x\in S^{01}_{\rm int}$.
Assume either condition {\rm (A)} or {\rm (B)} as follows.
\begin{itemize}
	\item[{\rm (A)}] $\eps=0$ and $\sum_{k} a(k)k^2<\infty$.
	\item[{\rm (B)}] $\eps\in(0,1)$ and $\sum_{k<0}a(k)k^2+\sum_{k>0} a(k)k<\infty$.
\end{itemize}
Then for any $t\geq0$,
\be\label{nonexp}
\E[h^\eps(X^\eps_0)]+\int_{0}^{t}\E\big[G^\eps h^\eps(X^\eps_s)\big]\di s\geq0,
\ee
where $h^\eps$ is the weighted number of inversions given in \eqref{hdef2}.
\el
\bl[Interface growth]\label{L:intgr}
Let $\sum_ka(k)|k|<\infty$, let $a$ be irreducible, and let $\eps\in[0,1)$.
Assume that interface tightness for $(X^\eps_t)_{t\geq0}$ in the sense of
Definition~\ref{D:int} does not hold on $S^{01}_{\rm int}$. Then the process
started in any initial state $X^\eps_0=x\in S^{01}_{\rm int}$ satisfies
\be\label{intgr}
\lim\limits_{T\to\infty}\frac{1}{T}\int_{0}^{T}\P\big[I_k(X^\eps_t)<N\big]\di t=0
\qquad(k>0,\ N<\infty),
\ee
where $I_k$ is given in (\ref{Idef}). The same statement holds with
$S^{01}_{\rm int}$ replaced by $S^{10}_{\rm int}$.
\el
With the lemmas above, we are ready to prove Theorem~\ref{T:tight}.\med

\bpro[~of Theorem~\ref{T:tight}]By symmetry, it suffices to consider the case with state space $S^{01}_{\rm int}$.
Let the biased voter model $(X^\eps_t)_{t\geq0}$ start from the heaviside state $x_0$ as in (\ref{heavy}).
Since $a$ is irreducible, for any constant $C$, there exist some $k_0\in\Z$ and $N\geq1$ such that
\be
C< \ha Na(k_0).
\ee
Let $C=\ha\sum_{k}a(k)k^2$ if $a$ has finite second moment, and $C=\ha\sum_{k<0}a(k)k^2+\eps^{-1}\sum_{k>0}a(k)k$ if $\eps>0$ and the moment condition
is relaxed to $\sum_{k<0}a(k)k^2<\infty$ and
$\sum_{k>0}a(k)k<\infty$.
Then for all $T>0$,
\bc\label{contra}
\dis 0 &\leq&\dis \frac{1}{T}\dis \int_{0}^{T}\E\big[G^\eps h^\eps(X^\eps_t)\big]\di t\\ [5pt]
&\leq&\dis C-\frac{1}{2T}\sum_{k}a(k)\int_{0}^{T}\E\big[I_k(X^\eps_t)\big]\di t
\leq C-\frac{Na(k_0)}{2T}\int_{0}^{T}\P\big[I_{k_0}(X^\eps_t)\geq N\big]\di t,
\ec
where in the first inequality we used Lemma~\ref{L:nonexp} and noted that
$\E[h^{\eps}(x_0)]=0$, and in the second inequality we used
Lemma~\ref{L:Geps}, in particular, the expression (\ref{Gh}) of $G^\eps
h^\eps$ and the inequalities $R^\eps_k\geq0$ and (\ref{Ghupp2}).
But on the other hand, if interface tightness did not hold, then by
Lemma~\ref{L:intgr},
\be
\lim\limits_{T\to\infty}\Big\{C-\frac{Na(k_0)}{2T} \int_{0}^{T}\P\big[I_{k_0}(X^\eps_t\geq N)\big]\di t\Big\}
= C-\ha N a(k_0)<0,
\ee
which contradicts (\ref{contra}).
Thus interface tightness must hold for the biased voter model.
\epro

\subsection{Proof of Lemma~\ref{L:Geps}}\label{S:Gh}

The proof is completed via a long calculation.
We first change some expressions into nice forms for later calculations.
Recall from (\ref{hdef2}) that for $\eps\in[0,1)$ and $x\in S^{01}_{\rm int}$,
\be
h^\eps(x)=\sum_{i}1_{\{x(i)=0\}}\sum_{j=-\infty}^{i-1}(1-\eps)^{\sum_{n<j}x(n)}1_{\{x(j)=1\}}.
\ee
Since the sum $\sum_{j=0}^{i-1}(1-\eps)^{j}$ is $\eps^{-1}\big(1-(1-\eps)^{i}\big)$ when $\eps>0$, we can rewrite $h^\eps(x)$ as
\be\label{hdef}
\ba{l}
\dis h^\eps(x)=\left\{\ba{ll}
\dis \sum_{i>j}\big(1-x(i)\big)x(j) & \dis\quad\dis\eps=0,\\ [5pt]
\dis\eps^{-1}\sum_i\big(1-x(i)\big)\big(1-(1-\eps)^{\sum_{j<i}x(j)}\big) & \dis\quad\dis\eps>0.
\ea\right.
\ec
For each $i\in\Z$, define functions
\be\label{fgdef}
\ba{l}
\dis f^{\eps}_i(x):=\left\{\ba{ll}
\dis\sum_{j<i}x(j)&\dis\mbox{if }\eps=0,\\[5pt]
\dis\eps^{-1}\big(1-(1-\eps)^{\sum_{j<i}x(j)}\big)\quad&\dis\mbox{if }\eps>0,
\ea\right. \\ [20pt]
\dis g_i(x):=1-x(i)=1_{\{x(i)=0\}},
\ec
and therefore
\be\label{hfg}
h^\eps=\sum_{i}f^\eps_ig_i.
\ee
We also rewrite the generator (\ref{bvmgen}) of the biased voter model as
\be\label{Gksum}
G^\eps=\sum_{k\neq 0}a(k)G^\eps_k,
\ee
where $G^\eps_k$ denotes the generator
\bc
\dis G^\eps_kf(x)
&:=&\dis\sum_n1_{\{x(n-k,n)=10\}}\big\{f(x+e_n)-f(x)\big\}\\[5pt]
&&+\dis(1-\eps)\sum_n1_{\{x(n-k,n)=01\}}\big\{f(x-e_n)-f(x)\big\},
\ec
which only describes infections $0\mapsto1$ and $1\mapsto0$ over distance $k$.

We start the calculations by recalling the following useful fact.
Let $X$ be a Markov process with countable state space $S$ and generator of the form
\be
Gf(x)=\sum_yr(x,y)\big\{f(y)-f(x)\big\},
\ee
where $r(x,y)$ is the rate of jumps from configuration $x$ to $y$.
Then for two real functions $f$ and $g$, by a direct calculation we have
\be\label{Gam}
G(fg)=fGg+gGf+\Ga(f,g),
\ee
where
\be
\Ga(f,g):=\sum_yr(x,y)\big\{f(y)-f(x)\big\}\big\{g(y)-g(x)\big\},
\ee
as long as all the terms involved are absolutely summable.

To find $G^\eps h^\eps$ for the biased voter model, applying formula (\ref{Gksum}) we can first calculate $G^\eps_k h^\eps$ and then sum over $k$.
By (\ref{hfg}) and (\ref{Gam}), we have
\bc\label{dcomp}
\dis G^\eps_kh^\eps&=&\dis\sum_iG^\eps_k(f^\eps_ig_i)
=\sum_i\big\{f^\eps_iG^\eps_kg_i+g_iG^\eps_kf^\eps_i-\Ga^\eps_k(f^\eps_i,g_i)\big\}\\[5pt]
&=&\dis\sum_i\big\{f^\eps_iG^\eps_kg_i+g_iG^\eps_kf^\eps_i\big\},
\ec
where in the last step we used that $\Ga^\eps_k(f^\eps_i,g_i)=0$, since any transition either changes the value of $x(i)$, in which case $f^\eps_i$ does not change, or the transition does not change the value of $x(i)$, in which case $g_i$ does not change.

We will prove that
\be\label{Gkh}
G^\eps_k h^\eps(x)
=\ha\big(k^2-I_k(x)\big)-\eps R^\eps_k(x).
\ee
Formula (\ref{Gh}) follows from this by summing over $k$ in $\Z$ with
weights given by $a$. We distinguish the calculation of
$G^\eps_kh^\eps$ into two cases, namely $k>0$ and $k<0$.
\medskip

\noindent {\bf Case $k>0$.} To ease notation, let us define a function $J_\eps:\N\to\R$ by
\be
J_\eps(n):=\left\{ \ba{ll}
n &\mbox{if }\eps=0,\\ [5pt]
\eps^{-1}\big(1-(1-\eps)^n\big)&\mbox{if }\eps>0,
\ea\right.
\ee
and thus for all $\eps\in[0,1)$,
\bc\label{lincr}
&\dis J_\eps(n+1)-J_\eps(n)=(1-\eps)^n,\\ [5pt]
&\dis f^\eps_i(x)=J_\eps\Big(\sum_{j<i}x(j)\Big).
\ec
Recall that $g_i(x)=1_{\{x(i)=0\}}$.
We have
\bc\label{GfGg}
\dis G^\eps_kf^\eps_i(x)
&=&\dis(1-\eps)\sum_{n<i}1_{\{x(n-k,n)=01\}}
\big\{J_\eps\big(\sum_{j<i}x(j)-1\big)-J_\eps\big(\sum_{j<i}x(j)\big)\big\}\\[5pt]
&&\dis+\sum_{n<i}1_{\{x(n-k,n)=10\}}
\big\{J_\eps\big(\sum_{j<i}x(j)+1\big)-J_\eps\big(\sum_{j<i}x(j)\big)\big\}\\[5pt]
&=&\dis-\sum_{n<i}1_{\{x(n-k,n)=01\}}(1-\eps)(1-\eps)^{\sum_{j<i}x(j)-1}\\[5pt]
&&\dis+\sum_{n<i}1_{\{x(n-k,n)=10\}}(1-\eps)^{\sum_{j<i}x(j)}\\[5pt]
&=&\dis(1-\eps)^{\sum_{j<i}x(j)}
\sum_{n<i}\big\{1_{\{x(n-k,n)=10\}}-1_{\{x(n-k,n)=01\}}\big\}\\[5pt]
&=&\dis(1-\eps)^{\sum_{j<i}x(j)}\sum_{n<i}\big\{1_{\{x(n-k)=1\}}-1_{\{x(n)=1\}} \big\}\\ [5pt]
&=&\dis-(1-\eps)^{\sum_{j<i}x(j)}\sum_{n=1}^k1_{\{x(i-n)=1\}},\\[5pt]
\dis G^\eps_kg_i(x)&=&\dis(1-\eps)1_{\{x(i-k,i)=01\}}-1_{\{x(i-k,i)=10\}}\\[5pt]
&=&\dis\big(1_{\{x(i-k)=0\}}-1_{\{x(i)=0\}}\big)-\eps1_{\{x(i-k,i)=01\}},
\ec
In the calculations above we used the equality
\be\label{reexp}
1_{\{x(i,j)=10\}}-1_{\{x(i,j)=01\}}=1_{\{x(i)=1\}}-1_{\{x(j)=1\}}=-1_{\{x(i)=0\}}+1_{\{x(j)=0\}},
\ee
which will be used repeatedly.
Substituting (\ref{GfGg}) into (\ref{dcomp}) leads to
\bc\label{sofar}
G^\eps_k h^\eps(x)
&=&\dis\sum_if^\eps_i(x)\big(1_{\{x(i-k)=0\}}-1_{\{x(i)=0\}}\big)-\eps\sum_if^\eps_i(x)1_{\{x(i-k,i)=01\}}\\[5pt]
&&\dis-\sum_i(1-\eps)^{\sum_{j<i}x(j)}\sum_{n=1}^k1_{\{x(i-n,i)=10\}}.
\ec
Note that in general if $f,g$ are functions such that $f_ig_i\to0$ as
$i\to\pm\infty$, then because
\be
f_ig_i-f_{i-n}g_{i-n}=\big(f_i-f_{i-n}\big)g_{i-n}+f_i\big(g_i-g_{i-n}\big),
\ee
one has the summation by parts formula
\be\label{sbp}
\sum_i\big(f_{i+n}-f_i\big)g_i=-\sum_if_i\big(g_i-g_{i-n}\big).
\ee
Applying this to $f_i(x)=f^\eps_i(x)$ and $g_i(x)=1_{\{x(i)=0\}}$, and using that $f^\eps_i(x)\to0$ as
$i\to-\infty$ and $1_{\{x(i)=0\}}\to0$ as $i\to\infty$, we have
\be\label{sofar0}
-\sum_if^\eps_i(x)\big(1_{\{x(i)=0\}}-1_{\{x(i-k)=0\}}\big)=\sum_i\big(f^\eps_{i+k}(x)-f^\eps_i(x)\big)1_{\{x(i)=0\}}.
\ee
We moreover observe that
by (\ref{lincr})
\be
J_\eps(n+k)-J_\eps(n)=\sum_{m=0}^{k-1}\big(J_\eps(n+m+1)-J_\eps(n+m)\big)=\sum_{m=0}^{k-1}(1-\eps)^{n+m},
\ee
and hence
\bc\label{sofar1}
\dis f^\eps_{i+k}(x)-f^\eps_i(x)
&=&\dis J_\eps\big(\sum_{j<i+k}x(j)\big)-J_\eps\big(\sum_{j<i}x(j)\big)\\[5pt]
&=&\dis\sum_{n=0}^{k-1}1_{\{x(i+n)=1\}}(1-\eps)^{\sum_{j<i+n}x(j)}.
\ec
Inserting the identities (\ref{sofar0}) and (\ref{sofar1}) into (\ref{sofar}) gives
\bc\label{sofar2}
G^\eps_kh^\eps(x)
&=&\dis\sum_i\sum_{n=0}^{k-1}1_{\{x(i,i+n)=01\}}(1-\eps)^{\sum_{j<i+n}x(j)}\\[5pt]
&&\dis-\sum_i(1-\eps)^{\sum_{j<i}x(j)}\sum_{n=1}^k1_{\{x(i-n,i)=10\}}\\[5pt]
&&\dis-\eps\sum_if^\eps_i(x)1_{\{x(i-k,i)=01\}}.
\ec
Before further simplifying (\ref{sofar2}), define counting functions $I^{10}_k,I^{01}_k:S^{01}_{\rm int}\to\N$ for all $k\in\Z$ by
\be\label{Idef2}
I^{10}_k(x):=\big|\{i:x(i,i+k)=10\}\big|\qquad\mbox{and}\qquad
I^{01}_k(x):=\big|\{i:x(i,i+k)=01\}\big|,
\ee
with $I^{10}_0=I^{01}_0:=0$. By the definition (\ref{Idef}) for the number of $k$-boundaries, we have $I_k=I^{10}_k+I^{01}_k$.
Moreover, we observe that for $x\in S^{01}_{\rm int}$, $k>0$ and $i\in\Z$, along the subsequence $x(\ldots,i-k,i,i+k,\ldots)$, there is one more adjacent pair $(01)$ than $(10)$, and thereby for any $k>0$,
\be\label{Irel}
I^{01}_k(x)=I^{10}_k(x)+k\qquad\mbox{and}\qquad
I_k(x)=2I^{10}_k(x)+k.
\ee
Changing the summation order, replacing $i$ by $i-n$ in the first term
in the right-hand side of (\ref{sofar2}), this term becomes
\be
%\weg{\sum_i\sum_{n=0}^{k-1}1_{\{x(i,i+n)=01\}}(1-\eps)^{\sum_{j<i+n}x(j)}=}\
\sum_i\sum_{n=0}^{k-1}1_{\{x(i-n,i)=01\}}(1-\eps)^{\sum_{j<i}x(j)},
\ee
we can rewrite (\ref{sofar2}) as
\bc\label{pos}
G^\eps_kh^\eps(x)
&=&\dis\sum_i(1-\eps)^{\sum_{j<i}x(j)}\sum_{n=1}^{k-1}
\big\{1_{\{x(i-n,i)=01\}}-1_{\{x(i-n,i)=10\}}\big\}\\[5pt]
&&\dis-\sum_i(1-\eps)^{\sum_{j<i}x(j)}1_{\{x(i-k,i)=10\}}\\[5pt]
&&\dis-\eps\sum_if^\eps_i(x)1_{\{x(i-k,i)=01\}}\\[5pt]
&=&\dis\sum_i\sum_{n=1}^{k-1}\big\{1_{\{x(i-n,i)=01\}}-1_{\{x(i-n,i)=10\}}\big\}\\[5pt]
&&\dis-\sum_i\big\{1-(1-\eps)^{\sum_{j<i}x(j)}\big\}\sum_{n=1}^{k-1}
\big\{1_{\{x(i-n,i)=01\}}-1_{\{x(i-n,i)=10\}}\big\}\\[5pt]
&&\dis-\sum_i1_{\{x(i-k,i)=10\}}
+\sum_i\big\{1-(1-\eps)^{\sum_{j<i}x(j)}\big\}1_{\{x(i-k,i)=10\}}\\[5pt]
&&\dis-\eps\sum_if^\eps_i(x)1_{\{x(i-k,i)=01\}}\\[5pt]
&=&\dis\sum_{n=1}^{k-1}\big(I^{01}_n(x)-I^{10}_n(x)\big)-I^{10}_k(x)\\[5pt]
&&\dis-\eps\sum_if^\eps_i(x)\sum_{n=1}^{k-1}
\big\{1_{\{x(i-n,i)=01\}}-1_{\{x(i-n,i)=10\}}\big\}\\[5pt]
&&\dis+\eps\sum_if^\eps_i(x)1_{\{x(i-k,i)=10\}}
-\eps\sum_if^\eps_i(x)1_{\{x(i-k,i)=01\}}\\[5pt]
&=&\dis\sum_{n=1}^{k-1}n-\ha(I_k(x)-k)-\eps\sum_if^\eps_i(x)\sum_{n=1}^k
\big\{1_{\{x(i-n,i)=01\}}-1_{\{x(i-n,i)=10\}}\big\}\\ [15pt]
&=&\dis\ha\big(k^2-I_k(x)\big)-\eps\sum_{n=1}^k\sum_if^\eps_i(x)\big(1_{\{x(i-n)=0\}}-1_{\{x(i)=0\}}\big).
\ec
where in the third and fourth equalities we used (\ref{Idef2}) and (\ref{Irel}), respectively, and in the last equality we interchanged the order of summation and used (\ref{reexp}).
We then substitute (\ref{sofar0}) back into the sum in the last line of (\ref{pos}), and use (\ref{sofar1}) to obtain
\be\label{rem}
\ba{l}
\dis\sum_{n=1}^{k}\sum_if^\eps_{i}(x)\big(1_{\{x(i-n)=0\}}-1_{\{x(i)=0\}}\big)
=\sum_{n=1}^{k}\sum_i\big(f^\eps_{i+n}(x)-f^\eps_{i}(x)\big)1_{\{x(i)=0\}}\\ [5pt]
\qquad=\dis\sum_{i}\sum_{n=1}^{k}\sum_{m=0}^{n-1}(1-\eps)^{\sum_{j<i+m}x(j)}1_{\{x(i,i+m)=01\}}\\ [5pt]
\qquad=\dis\sum_{i}\sum_{m=0}^{k-1}\sum_{n=m+1}^{k}(1-\eps)^{\sum_{j<i}x(j)}1_{\{x(i-m,i)=01\}}\\ [5pt]
\qquad=\dis\sum_{i}\sum_{m=1}^{k-1}(1-\eps)^{\sum_{j<i}x(j)}(k-m)1_{\{x(i-m,i)=01\}},
\ec
which implies (\ref{Gkh}) for $k>0$.
\medskip

\noindent {\bf Case $k<0$.} Similarly, we calculate
\bc
\dis G^\eps_kf^\eps_i(x)
&=&\dis(1-\eps)\sum_{n<i}1_{\{x(n,n+|k|)=10\}}
\big\{J_\eps\big(\sum_{j<i}x(j)-1\big)-J_\eps\big(\sum_{j<i}x(j)\big)\big\}\\[5pt]
&&\dis+\sum_{n<i}1_{\{x(n,n+|k|)=01\}}
\big\{J_\eps\big(\sum_{j<i}x(j)+1\big)-J_\eps\big(\sum_{j<i}x(j)\big)\big\}\\[5pt]
&=&\dis-\sum_{n<i}1_{\{x(n,n+|k|)=10\}}(1-\eps)(1-\eps)^{\sum_{j<i}x(j)-1}\\[5pt]
&&\dis+\sum_{n<i}1_{\{x(n,n+|k|)=01\}}(1-\eps)^{\sum_{j<i}x(j)}\\[5pt]
&=&\dis(1-\eps)^{\sum_{j<i}x(j)}
\sum_{n<i}\big\{1_{\{x(n,n+|k|)=01\}}-1_{\{x(n,n+|k|)=10\}}\big\}\\[5pt]
&=&\dis(1-\eps)^{\sum_{j<i}x(j)}\sum_{n=0}^{|k|-1}1_{\{x(i+n)=1\}},\\ [5pt]
\dis G^\eps_kg_i(x)&=&\dis(1-\eps)1_{\{x(i,i+|k|)=10\}}-1_{\{x(i,i+|k|)=01\}}\\[5pt]
&=&\dis\big(1_{\{x(i+|k|)=0\}}-1_{\{x(i)=0\}}\big)-\eps1_{\{x(i,i+|k|)=10\}},
\ec
which gives
\bc\label{sofar3}
G^\eps_kh^\eps(x)
&=&\dis\sum_if^\eps_i(x)\big(1_{\{x(i+|k|)=0\}}-1_{\{x(i)=0\}}\big)-\eps\sum_if^\eps_i(x)1_{\{x(i,i+|k|)=10\}}\\[5pt]
&&\dis+\sum_i(1-\eps)^{\sum_{j<i}x(j)}\sum_{n=1}^{|k|-1}1_{\{x(i,i+n)=01\}}.
\ec
Since by summation by parts,
\bc
\dis\sum_{i}f^\eps_i(x)\big(1_{\{x(i+|k|)=0\}}-1_{\{x(i)=0\}}\big)
&=& \dis-\sum_{i}\big(f^\eps_{i}(x)-f^\eps_{i-|k|}(x)\big)1_{x(i)=0}\\ [5pt]
&=&\dis-\sum_i\sum_{n=1}^{|k|}(1-\eps)^{\sum_{j<i-n}x(j)}1_{\{x(i-n,i)=10\}}\\ [5pt]
&=&\dis-\sum_{i}(1-\eps)^{\sum_{j<i}x(j)} \sum_{n=1}^{|k|}1_{\{x(i,i+n)=10\}},
\ec
combining the first and last terms on the right-hand side of (\ref{sofar3}), we can rewrite
\bc\label{negGh}
\dis G^\eps_kh^\eps(x)&=&
\dis\sum_{i} \sum_{n=1}^{|k|-1}\big\{1_{\{x(i,i+n)=01\}}-1_{\{x(i,i+n)=10\}}\big\}\\ [5pt]
&&\dis-\sum_{i}\big\{1-(1-\eps)^{\sum_{j<i}x(j)}\big\} \sum_{n=1}^{|k|-1}\big\{1_{\{x(i,i+n)=01\}}-1_{\{x(i,i+n)=10\}}\big\}\\ [5pt]
&&\dis-\sum_i 1_{\{x(i,i+|k|)=10\}}
\dis+\sum_{i}\big\{1-(1-\eps)^{\sum_{j<i}x(j)}\big\}
1_{\{x(i,i+|k|)=10\}} \\[5pt]
&&\dis-\eps\sum_if^\eps_i(x)1_{\{x(i,i+|k|)=10\}} \\ [5pt]
&=&\dis\sum_{n=1}^{|k|-1}\big(I^{01}_n(x)-I^{10}_n(x)\big)-I^{10}_{|k|}(x)\\ [5pt]
&&\dis-\eps\sum_if^\eps_i(x)\sum_{n=1}^{|k|-1}\big(1_{\{x(i)=0\}}-1_{\{x(i+n)=0\}}\big),\\ [5pt]
&=&\dis \ha \big(k^2-I_{|k|}(x)\big)-\eps R^\eps_k(x),
\ec
where in the last equality we have used the summation by parts formula (\ref{sbp}) and then (\ref{sofar1}) as follows
\be
\ba{l}
\dis\sum_{n=1}^{|k|-1}\sum_if^\eps_{i}(x)\big(1_{\{x(i)=0\}}-1_{\{x(i+n)=0\}}\big)
=\sum_{n=1}^{|k|-1}\sum_i\big(f^\eps_{i}(x)-f^\eps_{i-n}(x)\big)1_{\{x(i)=0\}}\\ [5pt]
\qquad=\dis\sum_{i}\sum_{n=1}^{|k|-1}\sum_{m=1}^{n}(1-\eps)^{\sum_{j<i-m}x(j)}1_{\{x(i-m,i)=01\}}\\ [5pt]
\qquad=\dis\sum_{i}\sum_{m=1}^{|k|-1}\sum_{n=m}^{|k|-1}(1-\eps)^{\sum_{j<i}x(j)}1_{\{x(i,i+m)=01\}}\\ [5pt]
\qquad=\dis\sum_{i}\sum_{m=1}^{|k|-1}(1-\eps)^{\sum_{j<i}x(j)}(|k|-m)1_{\{x(i,i+m)=01\}}.
\ec
Since $I_{-k}(x)=I_k(x)$ according to the definition~(\ref{Idef}), by
(\ref{negGh}), we see that (\ref{Gkh}) holds also for $k<0$.

In order to obtain the inequality (\ref{Ghupp2}) when $\eps\in(0,1)$, let $i_0:=\inf\{i\in\Z:x(i)=1\}$ and inductively let $i_n:=\inf\{i>i_{n-1}\in\Z:x(i)=1\}$, i.e., $i_0,i_1,\ldots$ are the positions of the first, second etc.\ $1$, coming from the left.
Thus, by counting from the left to right, for $k>0$,
\bc\label{Rlbb2}
\dis R^\eps_k(x)&=&\dis\sum_{n=0}^{\infty}(1-\eps)^{n}\sum_{m=1}^{k-1}(k-m)1_{\{x(i_n-m)=0\}}\\ [5pt]
&=&\dis\sum_{n=0}^{\infty}(1-\eps)^{n}\Big(\sum_{m=1}^{k-1}(k-m)-\sum_{m=1}^{k-1}(k-m)1_{\{x(i_n-m)=1\}}\Big)\\ [5pt]
&\geq&\dis\sum_{n=0}^{\infty}(1-\eps)^{n}\Big(\ha k(k-1)-kn\Big)\\ [5pt]
&=&\dis\ha\eps^{-1}k(k-1)-\eps^{-2}(1-\eps)k,
\ec
where in the inequality we bounded $\sum_{m=1}^{k-1}(k-m)1_{\{x(i_n-m)=1\}}$
by $k\sum_{m=1}^{k-1}1_{\{x(i_n-m)=1\}}$ and then used that there are at most
$n$ $1$'s on the left of site $i_n$, and in the last equality we
used the identity $\sum_{n=0}^{\infty}(1-\eps)^nn=\eps^{-2}(1-\eps)$.
Inserting \eqref{Rlbb2} into (\ref{Gkh}), we obtain, for $\eps>0$,
\be%\label{Gkhupp}
G^\eps_k h^\eps(x)
\leq\ha\big(k^2-I_k(x)\big)-\big(\ha k^2-\ha k -\eps^{-1}k+k \big)
\leq \eps^{-1}k-\ha I_k(x) \quad(k>0).
\ee
Using this and combining it for $k<0$ with the more elementary estimate
  $R^\eps_k(x)\geq0$ in (\ref{Gkh}), and summing over $k$, we then obtain
  (\ref{Ghupp2}).

\subsection{Proof of Lemma~\ref{L:nonexp}}\label{S:nonexp}

Though Lemma~\ref{L:nonexp} under condition~(A) is exactly Lemma~3 of \cite{SS08b}, we cannot follow the proof there to show our result under condition~(B).
More precisely, the estimate (3.24) in \cite{SS08b} cannot be used in case (B) due to the loss of finite second moment.
Instead, our proof uses different estimates, see Lemmas~\ref{L:Ibdd} and \ref{L:EIeps10} below, which turn out to work for the lemma under condition~(A) as well.

Let us recall from (\ref{Idef2}) that $I^{10}_k(x)=\big|\{i:x(i,i+k)=10\}\big|$.
\bl[Bound on number of inversions]\label{L:Ibdd}
Let $x\in S^{01}_{\rm int}$ and $I_n^{10}(x)$ be as in \eqref{Idef}.
Then for any $0\leq n<m$,
\be\label{Ibdd}
|I^{10}_m(x)-I^{10}_n(x)|\leq(m-n)I^{10}_1(x).
\ee
\el
\bpro
Suppose the interface of $x$ consists of $K$ blocks of consecutive 1's and $K$ blocks of consecutive 0's as follows.
\[
\cdots 0000000000
\underbrace{1111\cdots1111}_{\mbox{\scriptsize$1$st~1~block}}
\overbrace{0000\cdots0000}^{\mbox{\scriptsize$1$st~0~block}}
\cdots\cdots\cdots
\underbrace{1111\cdots1111}_{\mbox{\scriptsize$K$-th~1~block}}
\overbrace{0000\cdots0000}^{\mbox{\scriptsize$K$-th~0~block}}
1111111111\cdots
\]
It is straightforward to see that $I^{10}_1(x)=K$.
Suppose that the $k$-th block of consecutive 1's is from site $i_k$ to site $j_k$.
Then
\be
I_n^{10}(x)=\sum_{k=1}^{K}\sum_{s=i_k}^{j_k}1_{\{x(s+n)=0\}},
\ee
and therefore for $0<n<m$,
\be\label{Idiff2}
\big|I^{10}_m(x)-I^{10}_n(x)\big|\leq
\sum_{k=1}^{K}\big|\sum_{s=i_k}^{j_k}1_{\{x(m+s)=0\}}- \sum_{s=i_k}^{j_k}1_{\{x(n+s)=0\}}\big|.
\ee
To further bound the right hand side, note that for any $a,b_1,b_2,c\in\Z$ with $a<b_1,b_2<c$,
\be
\sum_{s=a}^{c}1_{\{x(s)=0\}}=\sum_{s=a}^{b_1-1}1_{\{x(s)=0\}}+\sum_{s=b_1}^{c}1_{\{x(s)=0\}}
=\sum_{s=a}^{b_2}1_{\{x(s)=0\}}+\sum_{s=b_2+1}^{c}1_{\{x(s)=0\}}.
\ee
Therefore
\be\label{rearr}
\sum_{s=b_1}^{c}1_{\{x(s)=0\}}-\sum_{s=a}^{b_2}1_{\{x(s)=0\}}
=\sum_{s=b_2+1}^{c}1_{\{x(s)=0\}}-\sum_{s=a}^{b_1-1}1_{\{x(s)=0\}}.
\ee
Applying this with $a=n+i_k,b_1=m+i_k,b_2=n+j_k$ and $c=m+j_k$ to (\ref{Idiff2}), we obtain

\bc\label{Idiff}
\big|I^{10}_m(x)-I^{10}_n(x)\big|&=&\dis\sum_{k=1}^{K}\big|\sum_{s=n+1}^{m}1_{\{x(j_k+s)=0\}} -\sum_{s=n}^{m-1}1_{\{x(i_k+s)=0\}}\big|\\ [5pt]
&\leq&\dis\sum_{k=1}^{K}\sum_{s=n+1}^{m}\big|1_{\{x(j_k+s)=0\}}
-1_{\{x(i_k+s-1)=0\}}\big|\\ [5pt]
&\leq&\dis\sum_{k=1}^{K}(m-n)=(m-n)I_1^{10}(x).
\ec
In particular, when $n=1<m$, (\ref{Idiff}) implies that
\be
|I^{10}_m(x)-I^{10}_1(x)|\leq(m-1)I^{10}_1(x),
\ee
which results in
\be
I^{10}_m(x)\leq mI^{10}_1(x).
\ee
The last inequality is nothing but (\ref{Ibdd}) for the case of $n=0<m$, and thus the proof is complete.
\epro

\medskip

\noi The following lemma  bounds uniformly the expected number of
$1$-boundaries $I_1$ from (\ref{Idef}).

%\bl\label{L:EI10}
%Let $(X^0_t)_{t\geq0}$ be the unbiased voter model starting from the heaviside state $x_0$ as in \eqref{heavy}.
%Then there exists a constant $C$ such that for all $k\in\Z$ and $t\geq0$,
%\be\label{EI1}
%\E\big[I^{10}_{k}(X^0_t)\big]\leq C|k|.
%\ee
%\el
%\bpro
%Let $B_{i,k,t}$ be the event that $X^0_t(i,i+k)=10$, that is, at time $t$ there is an inversion %pair $(10)$ at sites $i$ and $i+k$.
%Let $(S^{i}_s)_{s\geq0}$ be a random walk on $\Z$ with $S_0=i$.
%Then by the well-known duality between the voter model and coalescing random walks (see \cite[Chapter ~III~Section~4]{L85}), $B_{i,k,t}$ occurs if and only if there is a pair of backward coalescing random walks $(S^{i}_s)_{s\geq0}$ and $(S^{i+k}_s)_{s\geq0}$ with the same kernel $\tilde{a}$ where $\tilde{a}(k)=a(-k)$, such that $S^{i+k}_t<0\leq S^{i}_t$.
%Therefore,
%\bc\label{voterI}
%\dis\E\big[I^{10}_k(X_t)\big]&=&\dis\sum_{i}\E\big[1_{B_{i,k,t}}\big]=\dis\sum_{i}\P\big(S^{i+k}_t<0\leq S^{i}_t\big)\\ [5pt]
%&=&\dis\sum_{i}\sum_{j>0}\P\big(S^{k}_t=-i-j<-i\leq S^{0}_t\big)\\ [5pt]
%&=&\dis\sum_{j>0} \P\big(S^{0}_t-S^{k}_t\geq j\big)\\ [5pt]
%&\leq&\dis\E\big|S^{0}_t-S^{k}_t\big|,
%\ec
%where in the third equality we translated the system so that the starting position of the second random walk $S^{i}$ is the origin.
%By Lemma~2 of \cite{CD95}, for the pair of coalescing random walks $(S^0,S^k)$, there exists a constant $C$ such that for any $k\in\Z$ and $t\geq0$,
%\be
%\E\big|S^{0}_t-S^{k}_t\big|\leq C|k|.
%\ee
%Hence the inequality \eqref{EI1} follows.
%\epro

\bl[Bound on $1$-boundaries]\label{L:EIeps10}
Let $\sum_{k}a(k)|k|<\infty$.
Let $(X^\eps_t)_{t\geq0}$ be a biased voter model starting from a fixed configuration $x\in S^{01}_{\rm int}$. Then
%Then there exist constants $C_1,C_2$ such that
\be\label{EIeps1}
\sup_{\eps\in[0,1)}\E\big[I_{1}(X^\eps_t)\big]\leq I_1(x)e^{Ct},
\ee
where $C:=2\sum_{k\neq 0}a(k)|k-1|$.
\el
\br
Below in Lemma~\ref{L:1bdbd}, we also give a bound on the expected number
  of $1$-boundaries under the invariant law. Although the statements are similar, the proofs
 of Lemmas~\ref{L:EIeps10} and \ref{L:1bdbd} are completely different.
\er
\bpro
We will couple the process $\big(I_{1}(X^\eps_t)\big)_{t\geq0}$ to a branching
process $(Z_t)_{t\geq0}$ in such a way that $I_{1}(X^\eps_t)\leq Z_t$ for all
$t\geq0$. The left-hand side of (\ref{EIeps1}) can then be uniformly bounded
from above by the expectation of $Z_t$, which in turn can be bounded from above by
the right-hand side of (\ref{EIeps1}).  To see the coupling, note that
%whenever an infection of the biased voter model increases the number of
%$1$-boundaries, it must jump across at least one $1$-boundary and end at
%least a distance two from this $1$-boundary.
if site $i$ alters its type at time $t$ due to an infection
from site $j$ (without loss of generality, we may assume $i<j$),
and this increases the number of $1$-boundaries, then we must
have $X_{t-}(i-1)=X_{t-}(i)=X_{t-}(i+1)\neq X_{t-}(j)$.
Thus, there exists at least one $k\in(i+1,j]$ such that
$X_{t-}(k-1,k)$ is a $1$-boundary and the infection from $j$ to $i$ crosses
  this $1$-boundary and ends at least a distance two from this $1$-boundary (namely $k-i\geq2$).
For each $1$-boundary, the total rate
  of infections that cross it in this way is at most $\sum_{k\neq 0}a(k)|k-1|$. Since
a single infection increases the number of $1$-boundaries at most by 2,
  we can a.s.\ bound $I(X^\eps_t)$ from above by a branching process
  $(Z_t)_{t\geq0}$ started in $Z_0=I_1(x)$, where each particle produces two
  offspring with rate $\sum_{k\neq 0}a(k)|k-1|$, leading to the bound
  (\ref{EIeps1}).
%Hence, the total rate at which $I(X^\eps_t)$ increases is upper bounded by
%$2CI(X^\eps_t)$.  Therefore, letting $(Z_t)_{t\geq0}$ be the branching process
%with $Z_0=I_1(x)$ and increasing rate $2CZ_t$, there exists a natural coupling
%such that $I_{1}(X^\eps_t)\leq Z_t$ for all $t\geq0$.  This finishes the
%proof.
\epro

\medskip

\bpro[~of~Lemma~\ref{L:nonexp}]
By standard theory, if $f$ is bounded with bounded $G^\eps f$, then
\be\label{hmart}
M_t(f):=f(X^\eps_t)-\int_{0}^{t}G^\eps f(X^\eps_s)\di s \qquad (t\geq0)
\ee
is a martingale.
In particular, if the weighted number of inversions $h^\eps$ and $G^\eps h^\eps$ were bounded,  then
$M_t(h^\eps)$ would be a martingale and (\ref{nonexp}) would hold with equality.
However, $h^\eps$ is unbounded. But with a bit of extra work, we show next
that the inequality in (\ref{nonexp}) still holds. We use a truncation-approximation argument as in \cite{SS08b}.

Let $(X^{\eps,K}_t)_{t\geq0}$ be the process with the same initial state $X^{\eps,K}_0=X^\eps_0$ and
truncated kernel $a^K(k):=1_{\{|k|\leq K\}}a(k)$ $(k\in\Z)$.
%generator
%\bc\label{Ktrunc}
%G^{0,K} f(x)
%&=&\dis\sum_{i,j:|i-j|\leq K}a(j-i)1_{\{x(i,j)=10\}}\big\{f(x+e_j)-f(x)\big\}\\[5pt]
%&&\dis+\sum_{i,j:|i-j|\leq K}a(j-i)1_{\{x(i,j)=01\}}\big\{f(x-e_j)-f(x)\big\}.
%\ec
Recall that $L(x)$, defined in (\ref{intlen}),
%\be
%L(x)=\max\{i:x(i)=0 \}- \min\{i:x(i)=1 \}+1
%\ee
denotes the interface length of a configuration $x\in S^{01}_{\rm int}$.
Define stopping times
\be
\tau_{K,N}:=\inf\{t\geq0:L(X^{\eps,K}_t)>N \}\quad{\rm and}\quad
\tau_{N}:=\inf\{t\geq0:L(X^{\eps}_t)>N \}.
\ee
Let $G^{\eps,K}$ denote the generator from
	(\ref{bvmgen}) with the kernel $a$ replaced by $a^K$. For fixed $K$ and
$N$, since $L(X^{\eps,K}_{t})$ is bounded by $K+N$ for all
$t\leq\tau_{K,N}$, and $h^\eps(x)\leq h^0(x)$ where the latter is further bounded by $L(x)^2$, the interface length squared, we conclude that
\be
M^{K,N}_t:=h^\eps(X^{\eps,K}_{t\wedge\tau_{K,N}})-\int_{0}^{t\wedge\tau_{K,N}}G^{\eps,K}h^\eps(X^{\eps,K}_{s})\di s \quad (t\geq0)
\ee
is a martingale, which yields
\be\label{trunc}
0\leq \E[h^\eps(X^{\eps,K}_{t\wedge\tau_{K,N}})]=\E[h^\eps(X^{\eps,K}_0)]+\E\Big[\int_{0}^{t\wedge\tau_{K,N}}G^{\eps,K}h^\eps(X^{\eps,K}_{s})\di s\Big].
\ee
We will take limits as $N, K\to\infty$. Since $X^{\eps, K}_0=X^\eps_0$, the lemma would follow once we show
\be\label{Ghconv}
\lim\limits_{N\to\infty}\lim\limits_{K\to\infty}\E\Big[ \int_{0}^t 1_{\{s<\tau_{K,N}\}}G^{\eps,K}h^\eps(X^{\eps,K}_{s})\di s\Big]
=\E\Big[\int_{0}^{t}G^\eps h^\eps(X^\eps_{s})\di s\Big].
\ee

First, %{\yu \sout{ observe that } the first moment of $a$ being finite guarantees that}
we can couple $X^\eps$ and $X^{\eps,K}$ such that there exists a random $K_0$ so that if $K>K_0$, then $\tau_{K,N}=\tau_N$ and $X^{\eps,K}_t=X^\eps_t$ for all $t\leq\tau_N$. In particular, almost surely,
\be
\lim_{K\to\infty} \tau_{K,N} = \tau_N \quad \mbox{and} \quad \lim_{K\to\infty} X^{\eps, K}(s) 1_{\{s<\tau_{K,N}\}} = X^\eps(s)1_{\{s<\tau_N\}}.
\ee
Indeed, using the graphical representation in terms of resampling and selection
arrows described in Subsection~\ref{S:BN}, we can couple two processes with
kernel $a$ and $a^K$ by disallowing the use of arrows over a distance larger
than $K$ in the second process. Since $a$ has a finite first moment, the 
infections spread at a finite speed and the process
started in $S^{01}_{\rm int}$ stays in this space for all time, which in
particular implies that up to any finite time $T$, only finitely many arrows
have been used. We can then choose for $K_0$ the length of the longest arrow
that has been used up to time $\tau_N$, i.e., the largest distance over
which an infection has taken place in the process with kernel $a$.
%{\yu [To see the existence of such a coupling, note that $a(\cdot)$ having a finite first moment guarantees that infections spread at a finite speed. In particular, up to any finite time $T$, the space-time region that can trace through a sequence of infections back to the interface at time $0$ must be almost surely finite.
%It is then easy to choose $K_0$, which is just the size of the largest jump among all infections that can be traced back to the interface at time $0$.]}
%take a large interval $[-C,C]$ containing the interface, the first moment of $a$ being finite guarantees that before any finite time $T$ there are only finitely many infections from outside that can affect sites in $[C,-C]$. By choosing $T$ and $C$ large enough so that $\tau_N(\omega)\leq T$ and the interface is always contained in $[-C,C]$ up to $T$,

%\Swa{[I do not like the argument above. First of all, it does not really say
%    what the coupling is. Second, it is not completely clear to me what is
%    meant with the phrase `` the space-time region that can trace through a
%    sequence of infections back to the interface at time $0$''. I guess you
%    want to say something like ``trace its ancestory back'' but even then it
%    is not completey clear to me. I suggest the following argument.]}

Next we note that almost surely, $\lim_{N\to\infty}\tau_N=\infty$. Indeed, if
\be\label{asym}
a_{\rm s}(k):=\ha\big(a(k)+a(-k)\big)
\ee
denotes the symmetrization of $a(\cdot)$, then we can further couple $X^\eps$ and $X^{\eps,K}$ with a uni-directional random walk $S$ with increment rate $q(n):=\sum_{k= n}^{\infty}2a_{\rm s}(k)$ $(n\geq1)$ and $S_0=L(X^\eps_0)$, such
that $L(X^\eps_t)\leq S_t$ and $L(X^{\eps,K}_t)\leq S_t$ for all $t\geq0$.
%Since by the assumption on $a$, the kernel $q(\cdot)$ has \Jan{a} finite first moment $\sum_{n}q(n)n\leq 2\sum_{k}a(k)k^2<\infty$, i
It is then not hard to see that $\lim_{N\to\infty}\tau_N=\infty$ almost surely.

Now recall that, similar to the expression (\ref{Gh}) for $G^\eps h^\eps(x)$, we have
\be\label{Gh2}
 G^{\eps, K}h^\eps(X^{\eps,K}_{s}) 1_{\{s<\tau_{K,N}\}} = \sum_{k} 1_{\{s<\tau_{K,N}\}} 1_{\{|k|\leq K\}} a(k)\big(\ha k^2
-\eps R^\eps_k(X^{\eps,K}_{s}) - \ha I_{k}(X^{\eps,K}_{s})\big).
\ee
Since $\tau_{K,N}\to \tau_N$, $\tau_N\to\infty$, and $X^{\eps, K}(s) \to X^\eps(s)$, we see that
for each $s\in [0,t]$, $k\in\Z$, and almost surely, as first $K\to\infty$ and then $N\to\infty$, the summand above converges to
\be\label{Gh3}
 a(k)\big(\ha k^2
-\eps R^\eps_k(X^{\eps}_{s})\big)  - \ha a(k) I_{k}(X^{\eps}_{s}),
\ee
the sum of which over $k$ gives $G^\eps h^\eps(X^\eps_s)$. We will extend this pointwise convergence to the convergence of their integral with respect to $\E\big[\int_0^t \di s \sum_k \cdot \big]$ in \eqref{Ghconv}.

We first treat the contribution from $I_{k}(X^{\eps,K}_{s})$ in \eqref{Gh2}. Note that by (\ref{Irel}) and Lemma~\ref{L:Ibdd}, we have
\be
I_k(X^{\eps,K}_{s})=I^{10}_k(X^{\eps,K}_{s})+I^{01}_k(X^{\eps,K}_{s})\leq |k|\big(1+2I^{10}_1(X^{\eps,K}_{s})\big).
\ee
Recall that there exists a random $K_0$ such that if $K\geq K_0$, then $X^{\eps,K}_s=X^\eps_s$ for any $s\leq\tau_{K, N}$.
On the event $K_0>K$, we have $I_1^{10}(X^{\eps,K}_s)\leq L(X^{\eps,K}_s)+1\leq N+1$ for all $s<\tau_{K,N}$. Therefore
\be\label{Ibdd3}
\E\Big[1_{\{K_0>K\}}\int_{0}^{t \wedge \tau_{K,N}}\sum_{|k|\leq K}a(k)I_k(X^{\eps,K}_s)\di s\Big]
\leq P\big(K_0>K\big)\sum_{k}a(k)|k|\int_{0}^{t}\big(3+2N\big)\di s\,
\ee
which tends to zero in the limit of first $K\to\infty$ and then $N\to\infty$.

On the event $K_0\leq K$, because $X^{\eps,K}_s=X^\eps_s$ for $s<\tau_N$ and $\tau_{K,N}=\tau_N$, we have
\be\label{Ibdd2}
1_{\{K_0\leq K\}} 1_{\{s<\tau_{K,N}\}} 1_{\{|k|\leq K\}} a(k)I_k(X^{\eps,K}_s)
\leq a(k)|k| \big(1+2I^{10}_1(X^\eps_s)\big),
\ee
where the right hand side is integrable with respect to $\E\big[\int_0^t \di s \sum_k \cdot \big]$ by Lemma~\ref{L:EIeps10}, and the left hand side converges pointwise to $a(k)I_k(X^\eps_s)$ as first $K\to\infty$ and then $N\to\infty$. Therefore by dominated convergence, together with
\eqref{Ibdd3}, we obtain
\be\label{Gh4}
\lim\limits_{N\to\infty}\lim\limits_{K\to\infty}\E\Big[ \int_{0}^t 1_{\{s<\tau_{K,N}\}}\ha \sum_{|k|\leq K} a(k)I_k(X^{\eps, K}_s) \di s\Big]
=\E\Big[\int_{0}^{t} \ha \sum_k a(k) I_k(X^\eps_s) \di s\Big].
\ee

To treat the contribution from $\ha k^2-\eps R^\eps_k(X^{\eps,K}_{s})$, we note that
by the expression (\ref{Rdef}) for $R^\eps_k(x)$, it is easy to see that
\be
 \eps R_k^\eps(x) \leq \eps \sum_{i=1}^{\infty}(1-\eps)^{i-1} \sum_{n=1}^{|k|-1}n=\ha |k|(|k|-1),
\quad \mbox{hence}\quad  \ha k^2  -\eps R^\eps_k(X^{\eps,K}_{s})  \geq \ha |k|.
\ee
On the other hand, by the lower bound (\ref{Rlbb2}) on $R_k^\eps(x)$ when $k>0$ and the fact that $R_k^\eps(x)\geq 0$ when $k<0$, we have
\be\label{Rbdd2}
\ha k^2-\eps R^\eps_k(X^{\eps, K}_s) \leq  1_{\{k<0\}} \ha a(k)k^2+1_{\{k>0\}} \eps^{-1} a(k)k,
\ee
which is integrable with respect to $\E\big[\int_0^t \di s \sum_k \cdot \big]$ by either condition (A) or (B), while the left hand side converges pointwise to $\ha k^2-\eps R^\eps_k(X^{\eps}_s)$ as first $K\to\infty$ and then $N\to\infty$. Dominated convergence theorem can then be applied, which together with \eqref{Gh4}, implies \eqref{Ghconv}.
\epro

\br\label{R:nonexp}
If we assume the kernel $a(\cdot)$ has finite third moment, then we can prove that the process
$M_t(h^\eps)=h^\eps(X^\eps_t)-\int_{0}^{t}G^\eps h^\eps(X^\eps_s)\di s$ is a martingale.
To prove this, we only need to check the uniform integrability of
$\big(h^{\eps}(X^{\eps,K}_{t\wedge\tau_{K,N}})\big)$ in $K$ and $N$, because then we
can take the limit on the left-hand side of (\ref{trunc}) as well, and the equality in (\ref{trunc}) remains valid as $K\to\infty$ and then $N\to\infty$.
Recall that the uni-directional random walk $S$ has increment rate $q(n)=\sum_{k=n}^{\infty}2a_{\rm s}(n)$ whose second moment is now finite since
\be
\E\big[(S_t-S_0)^2\big]\leq t\sum_{n}q(n)n^2\leq t\sum_{k}a(k)|k|^3<\infty.
\ee
The uniform integrability thus follows from the fact that
\be
h^{\eps}(X^{\eps,K}_{t\wedge\tau_{K,N}})\leq
h^0(X^{\eps,K}_{t\wedge\tau_N})\leq L(X^{\eps,K}_{t\wedge\tau_{K,N}})^2
\leq S_t^2 \quad{\rm a.s.},
\ee
where in the third inequality we used $L(X^{\eps,K}_{t\wedge\tau_N})\leq S_{t\wedge\tau_N}\leq S_t$.

In particular, for the unbiased voter model $(X^0_t)_{t\geq0}$, the process $(M_t(h^0))_{t\geq0}$ is a martingale.
We claim, however, if we replace $h^0$ by $h_M$, the number of inversions within distance $M$, formally defined by
\be\label{htrun}
h_M(x):=\big|\{(j,i):0<i-j\leq M,~x(j,i)=10\}\big|,
\ee
then a finite second moment assumption would suffice to imply that $M_t(h_M)$ is a martingale, and therefore
\be\label{hMmart0}
\E\big[h_M(X^0_t)\big]-\E\big[h_M(X^0_0)\big]=\E\Big[\int_{0}^{t}G^0h_M(X^0_s) \di s\Big].
\ee
Indeed, since the inversion pairs must be inside the interface and each particle in the interface contributes to at most $2M$ pairs of inversions, for any $x\in S^{01}_{\rm int}$ we have
\be
h_M(x)\leq 2ML(x).
\ee
Thus the uniform integrability of
$\big(h_M(X^{0,K}_{t\wedge\tau_{K,N}})\big)$ follows from
\be
h_M(X^{0,K}_{t\wedge\tau_{K,N}})\leq2ML(X^{0,K}_{t\wedge\tau_{K,N}})\leq 2MS_t\quad{\rm a.s.}\qquad{\rm and}\qquad \E[S_t]<\infty,
\ee
which only requires $a$ to have finite second moment.
Therefore, in order to show that $M_t(h_M)$ is a martingale, it remains to check the uniform integrability of $\int_{0}^{t\wedge\tau_{K,N}}G^{0,K}h_M(X^{0,K}_s) \di s$.
Using the expression of $G^{0,K}h_M$ in (\ref{Ghhat0}) below, one only needs to show the uniform integrability of
\be\label{Itrun}
\int_{0}^{t\wedge\tau_{K,N}}\left(\sum_{n=1}^{\infty}A^K(n)I^{10}_{M+n}(X^{0,K}_s)-\sum_{n=1}^{\infty}A^K(n)I^{10}_{M-1-n}(X^{0,K}_s)\right)\di s \qquad(K,N\geq1)
\ee
where $A^K(n)=\sum_{k=n}^{\infty}\big(a^K(k)+a^K(-k)\big)$.
Estimating $I^{10}_{M+n}\leq I_{M+n}\leq(M+n)I_1$ and $I^{10}_{M-1-n}\leq I_{M-1-n}\leq|M-1-n|I_1$ and using Lemma~\ref{L:EIeps10}, one gets that the first moment of $A(\cdot)$, or equivalently, the second moment of $a(\cdot)$, being finite guarantees the uniform integrability of the terms in (\ref{Itrun}).
Thus the equality (\ref{hMmart0}) has been proved, which we state as the following lemma.
This result will be used in the proof of Proposition~\ref{P:equi} later on.
\er
\bl\label{L:hMmart}
Let the voter model $(X^0_t)_{t\geq0}$ start from a fixed configuration $X^0_0=x\in S^{01}_{\rm int}$, and let $h_M$ denote the number of inversions within distance $M$ as in \eqref{htrun}.
Assume that $\sum_{k} a(k)k^2<\infty$.
Then for any $t\geq0$,
\be\label{hMmart}
\E\big[h_M(X^0_t)\big]-\E\big[h_M(X^0_0)\big]=\E\Big[\int_{0}^{t}G^0h_M(X^0_s) \di s\Big].
\ee
\el

\subsection{Proof of Lemma~\ref{L:intgr}}\label{S:intgr}

We fix $\eps\in[0,1)$ throughout the proof. We only state the proof for
$S^{01}_{\rm int}$, the proof for $S^{10}_{\rm int}$ being the same.
We start by proving the statement for $k=1$. Let $(X^\eps_t)_{t\geq 0}$ be
started in an initial state $X^\eps_0=x\in S^{01}_{\rm int}$ and consider the
``boundary process'' $(Y_t)_{t\geq 0}$ defined as
\be
Y_t(i):=X^\eps_t(i+1)-X^\eps_t(i)\qquad(i\in\Z,\ t\geq 0).
\ee
The assumption $\sum_ka(k)|k|<\infty$ guarantees that a.s.\ $X^\eps_t\in
S^{01}_{\rm int}$ for all $t\geq 0$ and hence $(Y_t)_{t\geq 0}$ is a Markov
process in the space of all configurations $y\in\{-1,0,1\}^\Z$ such that
$\sum_i|y(i)|$ is odd (and finite) and $\sum_{i:i\leq j}y(i)\in\{0,1\}$ for
all $j\in\Z$. For any such configuration, we set $|y|:=\sum_i|y(i)|$.
Then $|Y_t|=I_1(X^\eps_t)$ $(t\geq 0)$.

In the special case that $\eps=0$, the process $(Z_t)_{t\geq 0}$ defined as
$Z_t(i):=|Y_t(i)|$ is also a Markov process, and in fact a cancellative spin
system. In this case, we can apply \cite[Proposition~13]{SS08a} to conclude that
\be\label{boundgrw}
\lim_{T\to\infty}\frac{1}{T}\int_0^T\P\big[|Y_t|\leq N\big]\di t=0
\qquad(N<\infty).
\ee
We describe this in words by saying that for each $N$, the process spends a
zero fraction of time in states $y$ with $|y|\leq N$. In the biased case
$\eps>0$, the process $(Z_t)_{t\geq 0}$ is no longer a Markov process, but we
claim that the proof of \cite[Proposition~13]{SS08a} can easily be adapted to
show that (\ref{boundgrw}) still holds. To demonstrate this, we go through the
main steps of that proof and show how to adapt them to our process $(Y_t)_{t\geq
  0}$.

The main ingredient in the proof of \cite[Proposition~13]{SS08a} is formula
(3.54) of that paper, which for our process must be reformulated as
\be\label{boundan}
\inf\big\{\P^y[|Y_t|=n]:|y|=n+2,\ y(i)\neq 0\neq y(j)
\mbox{ for some }i\neq j,\ |i-j|\leq L\big\}>0
\ee
for all $t>0$, $L\geq 1$, and $n=1,3,5,\ldots$. Let us call a site $i\in\Z$
such that $Y_t(i)\neq 0$ a boundary of $X^\eps_t$. Then (\ref{boundan}) says
that if $X^\eps_t$ contains $n+2$ boundaries of which two are at distance $\leq
L$ of each other, then there is a uniformly positive probability that after
time $t$ the number of boundaries has decreased by 2.

The assumption that interface tightness for $X^\eps_t$ does not hold on
$S^{01}_{\rm int}$ implies that (\ref{boundgrw}) holds for $N=1$. The proof of
(\ref{boundgrw}) now proceeds by induction on $N$. Imagine that
(\ref{boundgrw}) holds for $N$. Then it can be shown that (\ref{boundan})
implies that for each $L\geq 1$, the process spends a zero fraction of time in
states $y$ with $|y|=N+2$ which contain two boundaries at distance $\leq L$ of
each other. Now imagine that (\ref{boundgrw}) does not hold for $N+2$. Then
for each $L\geq 1$, the process must spend a positive fraction of time in
states $y$ with $|y|=N+2$ but which do not contain two boundaries at distance
$\leq L$ of each other.

If $L$ is large, then each boundary evolves for a long time as a process
started in a heaviside initial state, either of type $01$ or of type $10$,
without feeling the other boundaries, which are far away. By our assumption
that interface tightness on $S^{01}_{\rm int}$ does not hold, the boundaries
of type $01$ are unstable in the sense that they will soon split into three or
more boundaries and on sufficient long time scales spend most of their
  time being three or more boundaries, rather than one. With some care, it
can be shown that this implies that the process spends a zero fraction of time
in states $y$ with $|y|=N+2$, completing the induction step. This argument is
written down more carefully in \cite{SS08a}. Formula (3.65) of that paper has
to be slightly modified in our situation since we only know that the boundaries
of type $01$ are unstable. So instead of producing at least $3(N+2)$
boundaries with probability at least $(1-2\eps)^{N+2}$,
in our case, we produce at least $3(N+3)/2+(N+1)/2$
boundaries with probability at least $(1-2\eps)^{(N+3)/2}$, since of the
$N+2$ boundaries there are $(N+3)/2$ of type $01$ and $(N+1)/2$ of type $10$.

Translated back to the biased voter model $(X^\eps_t)_{t\geq 0}$ started
in an initial state $X^\eps_0=x\in S^{01}_{\rm int}$, formula (\ref{boundgrw})
says that
\be\label{onestep}
\lim\limits_{T\to\infty}\frac{1}{T}\int_{0}^{T}\P[I_1(X^\eps_t)<N]\di t=0
\qquad(N<\infty).
\ee
To deduce (\ref{intgr}) from (\ref{onestep}), it suffices to prove that for any $k,M\geq1$ and $s>0$,
\be\label{mulstep}
\lim\limits_{N\to\infty}~\!\sup\limits_{X^\eps_0:I_1(X^\eps_0)\geq N}\P[I_k(X^\eps_s)<M]=0.
\ee
For if (\ref{mulstep}) holds, then for any $s,\delta>0$, there exists $N$
large enough such that the supremum in (\ref{mulstep}) is less than $\delta$,
which, letting the process in (\ref{onestep}) evolve for some extra time $s$ implies that
\be
\limsup\limits_{T\to\infty}\frac{1}{T}\int_{0}^{T}\P[I_k(X^\eps_{t+s})<M]\di t< 2\delta \quad (M\geq1).
\ee
As $\delta$ is arbitrary, (\ref{intgr}) is obtained.

It remains to show (\ref{boundan}) and (\ref{mulstep}). Both of them can be
proved by directly constructing specific paths with positive probabilities, by
constructions very similar to those in the proof of \cite[Proposition~4]{SS08b}.

\section{Continuity of the invariant law}\label{S:nucont}

\subsection{Proof outline}\label{S:contout}

As a positive recurrent Markov chain on the countable state space
$\ov{S}^{01}_{\rm int}$, the biased voter model modulo translations has a
unique invariant law $\ov{\nu}_\eps$. This section is devoted to proving
Theorem~\ref{T:nucont}, namely the weak convergence
$\ov{\nu}_\eps\Rightarrow\ov{\nu}_0$ with respect to the discrete topology
on~$\ov{S}^{01}_{\rm int}$.

%\sout{
%Recall that $\ov{S}^{01}_{\rm int}$ is the set of equivalence classes of
%elements of $S^{01}_{\rm int}$ that are equal up to translations. It will be
%convenient to choose a representative from each equivalence class by shifting
%the leftmost one to the origin. Since each equivalence class $\ov x\in\ov
%S^{01}_{\rm int}$ contains a unique element $x$ in the set
%we can identify $\ov S^{01}_{\rm int}$ with $\hat S^{01}_{\rm int}$. Under
%this identification, $\ov{\nu}_\eps$ on $\ov{S}^{01}_{\rm int}$ uniquely
%determines a probability measure $\nu_\eps$ on $\{0,1\}^\Z$ that is supported
%on $\hat{S}^{01}_{\rm int}$. We let $X^\eps_\infty$ denote a random variable
%with law~$\nu_\eps$.
%}
Recall from Section \ref{S:intrnucon} that for each $\eps\in[0,1)$, by selecting
a representative configuration with the leftmost 1 at the origin, $\ov{\nu}_\eps$ uniquely
determines a probability measure $\nu_\eps$ on $\{0,1\}^\Z$ that is supported
on 
\be\label{hatS}
\hat{S}^{01}_{\rm int}:=\{x\in S^{01}_{\rm int}:x(i)=0
\mbox{ for all $i<0$ and }x(0)=1\}.
\ee
Let $X^\eps_\infty$ denote a random variable
with law~$\nu_\eps$.
If we could show tightness for the length of the interface $L(X^\eps_\infty)$
(see (\ref{intlen})) as $\eps$ varies, then it would be relatively straightforward to
prove weak continuity of the map $\eps\mapsto\nu_\eps$ with respect to the
discrete topology on $\hat S^{01}_{\rm int}$. Unfortunately, our proof of the existence of 
the equilibrium interface $X^\eps_\infty$ gives us little control on
$L(X^\eps_\infty)$. However, we are able to derive a uniform upper bound on the expected 
number of weighted $k$-boundaries (see (\ref{Idef})) in the equilibrium biased voter interface $X^\eps_\infty$. As a first
step, this is sufficient to prove weak continuity of the map $\eps\mapsto\nu_\eps$ 
with respect to the {\em product topology} on~$\hat S^{01}_{\rm int}$, which implies in particular
that $\nu_\eps\Rightarrow \nu_0$ as $\eps\downarrow 0$ under the product topology.

To boost this to
weak convergence with respect to the discrete topology on $\hat{S}^{01}_{\rm int}$, from which Theorem \ref{T:nucont} then follows, we argue by contradiction.
%We do not know if the map $\eps\mapsto\nu_\eps$ is continuous with respect to the discrete topology on $\hat S^{01}_{\rm int}$, but we can establish continuity at $\eps=0$ by a rather subtle argument.
Let $\eps_n\down 0$ and
assume that $\nu_{\eps_n}\Rightarrow\nu_0$ with respect to the product
topology on $\hat S^{01}_{\rm int}$, but not with respect to the discrete
topology. Then, with positive probability, $X^{\eps_n}_\infty$ must contain
boundaries that ``walk to infinity'' as $\eps_n\down 0$. We will rule out this scenario
by proving that the expectation of a weighted sum of the $k$-boundaries in the equilibrium biased
voter model interface cannot exceed that of the voter model (see Lemma~\ref{L:equibdd}). The latter can be computed from the equilibrium equation $\E[G^0h^0(X^0_\infty)]=0$ and equals $\frac{1}{2}\sig^2$ (see Proposition \ref{P:equi}).

The rest of this section is organized as follows. In
Subsection~\ref{S:bndbnd} we derive a uniform upper bound on the expected
number of weighted $k$-boundaries in the equilibrium biased voter interface.  
In Subsection~\ref{S:ConPT}, we
prove weak continuity of the map $\eps\mapsto\nu_\eps$ with
respect to the product topology on $\hat S^{01}_{\rm int}$. In
Subsection~\ref{S:equi}, we establish the equilibrium equation
$\E[G^0h^0(X^0_\infty)]=0$, which determines the expected weighted sum of $k$-boundaries
in the equilibrium voter interface. This is the most technical part due to the unboundedness of
$h^0$. Once this is done, however, the proof of Theorem~\ref{T:nucont} is quite short and is given in Subsection~\ref{S:ConS}.

\subsection{Bound on the number of boundaries}\label{S:bndbnd}

The following lemma is our basic tool to control the expected weighted sum of  of boundaries
in the equilibrium biased voter interface. We note that in Proposition~\ref{P:equi} below,
we will show that \eqref{equibdd} is in fact an equality when $\eps=0$. Thus the 
expected weighted sum of of boundaries of the equilibrium voter interface gives an upper bound for that of 
the biased voter model. This fact will be key to our proof of Theorem~\ref{T:nucont}.

\bl[Bound on $k$-boundaries]\label{L:equibdd}
For $\eps\in[0,1)$, let $X^\eps_\infty$ denote a random variable
with law~$\nu_\eps$ as in Subsection~\ref{S:contout}. Then
\be\label{equibdd}
\E\big[\sum_{k=1}^{\infty}a_{\rm s}(k)I_k(X^\eps_\infty)\big]\leq\ha\sig^2,
\ee
where $I_k$ is defined in \eqref{Idef}, $a_{\rm s}(k)=\ha\big(a(k)+a(-k)\big)$ and $\sig^2=\sum_{k\in\Z}a(k)k^2$.
\el
\bpro
Let the biased voter model $(X^\eps_t)_{t\geq0}$ start from the heaviside initial state $X^\eps_0=x_0$ as in \eqref{heavy}.
For $\eps\in(0,1)$, recall from \eqref{hdef2} that
\be
h^\eps(x)=\sum_{i>j}(1-\eps)^{\sum_{n<j}x(n)}1_{\{x(j,i)=10\}}.
\ee
Hence $h^\eps(X^\eps_0)\equiv0$.
By Lemma~\ref{L:Geps}~and Lemma~\ref{L:nonexp}, for any $t>0$,
\be\label{equibdd3}
\int_{0}^{t}\E\big[\sum_{k=1}^{\infty}a_{\rm s}(k)\big(k^2-I_k(X^\eps_s)\big)\big]\di s\geq
\int_{0}^{t}\E\big[G^\eps h^\eps(X^\eps_s)\big]\di s\geq 0,
\ee
or equivalently,
\be
\label{equibdd2}
\int_{0}^{t}\E\big[\sum_{k=1}^{\infty}a_{\rm s}(k)I_k(X^\eps_s)\big]\di s\leq
\int_{0}^{t}\E\big[\sum_{k=1}^{\infty}a_{\rm s}(k)k^2\big]\di s=\frac{1}{2}\sig^2 t.
\ee
Dividing both sides by $t$ and then letting $t\to\infty$, we arrive at
(\ref{equibdd}) by Fatou's lemma, since $\ov{X}^\eps_t$ converges weakly to
$\ov{X}^\eps_\infty$. 
%\Swa{[I think we still need to add a bit of detail since Fatou
%    usually refers to integrals for which the integrand converges
%    pointwise. Here the setting is different since we have weak convergence.]}
\epro

\bl[Bound on $1$-boundaries]
There\label{L:1bdbd} exists a constant $C<\infty$ such that
\be\label{1bdbd}
\sup_{\eps\in[0,1)}\E\big[I_1(X^\eps_\infty)\big]\leq C.
\ee
\el
\bpro
Fix $t>0$, and choose $k\geq 1$ such that $a_{\rm s}(k)>0$. For a biased voter
model started in any initial state $x$ with $x(i,i+1)=10$, using the
irreducibility of the kernel $a$, it is easy to see that there is a positive
probability $p$ that the 1 at position $i$ spreads to position $i-k+1$ at time $t$, while leaving
the 0 at position $i+1$ as it is. Since this event only requires the 1's to spread, this
probability can be bounded from below uniformly in the bias $\eps$. Therefore, if
$(X^\eps_t)_{t\geq 0}$ denotes the biased voter model started with initial law $\nu_\eps$, then
\be\label{I1bdd}
\E\big[I^{10}_1(X^\eps_0)\big]=\sum_{i}\E\big[1_{X^\eps_0(i,i+1)=10}\big]
\leq\sum_{i}\frac{1}{p}\E\big[1_{X^\eps_t(i-k+1,i+1)=10}\big]
=\frac{1}{p}\E\big[I^{10}_k(X^\eps_t)\big].
\ee
Since the law of $X^\eps_t$ modulo translations does not depend on $t$,
the claim \eqref{1bdbd} follows from Lemma~\ref{L:equibdd}.
\epro

\subsection{Continuity with respect to the product topology}\label{S:ConPT}

Throughout this subsection, we will consider the biased voter model $(X^\eps_t)_{t\geq0}$ with initial law $\nu_\eps$, which corresponds to the equilibrium interface with the leftmost $1$ shifted to the origin. We will prove the following theorem.

\bt[Continuity with respect to the product topology]
Assume\label{T:convpt} that the kernel $a(\cdot)$ is irreducible and satisfies
$\sum_ka(k)k^2<\infty$. Then the map $[0,1)\ni\eps\mapsto\nu_\eps$ is continuous with respect to weak
convergence of probability measures on $\{0,1\}^\Z$, equipped with the product topology. 
\et

To prepare for the proof of Theorem~\ref{T:convpt}, we need a few lemmas.
%Let $(X^\eps_t)_{t\geq 0}$ denote the biased voter model started in the initial state $\nu_\eps$. 
At time $t\geq0$, denote the position of the leftmost $1$ and the rightmost 0 %\sout{zero}
by $l^\eps_t$ and $r^\eps_t$, respectively. That is,
\be\label{left}
l^\eps_t:=\min\{i:X^\eps_t(i)=1\}
\quand
r^\eps_t:=\max\{i:X^\eps_t(i)=0\}\qquad(t\geq 0).
\ee
For $i\in\Z$, define a shift operator $\tet_i : \{0,1\}^\Z \to \{0, 1\}^\Z$ by
\be
\tet_i(x)(j):=x(i+j)\qquad(i,j\in\Z,\ x\in S^{01}_{\rm int}).
\ee
Then $\tet_{l^\eps_t}(X^\eps_t)$ has law $\nu_\eps$ for all $t\geq 0$.

Our strategy for proving Theorem~\ref{T:convpt} is as follows. Fix
$\eps_n,\eps^\ast\in[0,1)$ such that $\eps_n\to\eps^\ast$. Since 
$\{0,1\}^\Z$ is compact under the product topology, tightness comes for free. So by going to a subsequence
if necessary, we can assume that $\nu_{\eps_n}\Rightarrow\nu^\ast$ for some
probability measure $\nu^\ast$ on $\{x\in\{0,1\}^\Z:x(0)=1,\ x(i)=0\ \forall i<0\}$.
%$\{0,1\}^\Z$. By Lemma~\ref{L:1bdbd}, we see that 
%each subsequential limit $\nu^\ast$ is concentrated on $\{ x\in S^{01}_{\rm
%  int} \cup S^{00}_{\rm int}: x(0)=1\}$, see~\eqref{eq:S01}. 
Using convergence of the generators, general arguments tell us that
$(X^{\eps_n}_t)_{t\geq0}$ converges in distribution as a process to the voter model $(X^\ast_t)_{t\geq0}$
with initial law $\nu^\ast$. This is Lemma~\ref{L:bvtcon} below. Next, in
Lemma~\ref{L:ltight}, we show that for any fixed $t$, the family
$(l^\eps_t)_{\eps\in[0,1)}$ is tight. By going to a further subsequence if necessary, this 
implies that $\tet_{l^\eps_t}(X^\eps_t) \Rightarrow \theta_{l^*_t}(X^*_t)$, where $l^*$ is the position of the leftmost 1 in $X^*_t$. 
Since $\tet_{l^\eps_t}(X^\eps_t)$ has law $\nu_\eps$ and converges to $\nu^*$, it follows that 
$\nu^*$ is an invariant law for $X^\ast_\cdot$ seen from the leftmost one.
By Lemma~\ref{L:1bdbd}, we see that $\nu^\ast$ is concentrated on $S^{01}_{\rm
	int} \cup S^{00}_{\rm int}$, see~\eqref{eq:S01}. It is not hard to see that
$\nu^\ast$ gives zero measure to $S^{00}_{\rm int}$,
%In fact, it is not hard to see that $\nu^*$ concentrates on $S_{\rm int}^{01}$, 
and hence must equal $\nu_{\eps^\ast}$ by the uniqueness of the invariant law
for the biased voter interface. This establishes the continuity of $\eps\mapsto\nu_\eps$.

%\Swa{[The outline above is not clear to me. It is clear that we can get a
%subsequence that converges to a law $\nu^\ast$ on $S_{\rm
%  left}:=\{x\in\{0,1\}^\Z:x(i)=0\ \forall i<0,\ x(0)=1\}$. By soft arguments,
%we can prove that $\nu^\ast$ is invariant for the biased voter model as seen
%from the left-most one. If we can prove that the latter has only one invariant
%law on the space $S_{\rm left}$, then we are done. The outine above, however,
%seems to suggests we know a priori that $\nu^\ast$ is concentrated on
%$S^{01}_{\rm int}$ and use uniqueness of the invariant law for the process on
%$S^{01}_{\rm int}$ as seen from the leftmost one, which follows from interface
%tightness on $S^{01}_{\rm int}$, which has already been proved.

%When I read the proof, it seems to be OK. We use Lemma~\ref{L:1bdbd} to
%conclude that $\nu^\ast$ is concentrated on the union of $S^{01}_{\rm int}$
%and $S^{00}_{\rm int}$. If the restriction of $\nu^\ast$ to any of these
%spaces is nonzero, then normalizing it gives an invariant law of the process
%as seen from the leftmost one. But it is not hard to see that this process
%cannot have an invariant law on $S^{00}_{\rm int}$. So we obtain that
%$\nu^\ast$ is concentrated on $S^{01}_{\rm int}$ and we are done.]}

\bl[Convergence as a process]
Assume\label{L:bvtcon} $\eps_n\rightarrow\eps^\ast\in[0,1)$ such that
$\nu_{\eps_n}\Rightarrow\nu^\ast$ for some probability measure $\nu^\ast$ on
$\{0,1\}^\Z$, equipped with the product topology. Then
\be\label{bvtcon}
\P\big[(X^{\eps_n}_t)_{t\geq0}\in\,\cdot\,\big]\Asto{n}
\P\big[(X^\ast_t)_{t\geq0}\in\,\cdot\,\big],
\ee
where $X^{\eps_n}_0$ has law %\sout{$\nu^{\eps_n}$}
$\nu_{\eps_n}$, $(X^\ast_t)_{t\geq0}$ is the biased voter model with bias $\eps^\ast$ and
initial law $\nu^\ast$, and $\Rightarrow$ denotes weak convergence in the Skorohod space $D([0,\infty),\{0,1\}^\Z)$.
\el
\bpro
By \cite[Theorem~3.9]{L85}, for each $\eps\in[0,1)$, the generator $G^\eps$ in
(\ref{bvmgen}) is well-defined for any function $f\in\Di$ with
\be
\Di:=\Big\{f:\sum_{i\in\Z}\sup_{x\in\{0,1\}^\Z}
\big|f(x+e_i)-f(x)\big|<\infty\Big\},
\ee
and the closure of the generator $G^\eps$ with domain $\Di$ generates a Feller
semigroup. By \cite[Theorem~17.25]{Kal97}, to establish (\ref{bvtcon}), it
suffices to check that $\|G^{\eps_n}f-G^{\eps^\ast}f\|_\infty\to 0$ for all
$f\in\Di$, where $\|\,\cdot\,\|_\infty$ denotes the supremum norm.
This follows by writing
\bc
\dis\big|G^{\eps_n}f(x)-G^{\eps^\ast}f(x)\big|
&=&\dis|\eps_n-\eps^\ast|\cdot\Big|\sum_{k\neq0}a(k)\sum_{i}1_{x(i-k,i)=01}
\big\{f(x-e_i)-f(x)\big\}\Big|\\ [5pt]
&\leq&\dis|\eps_n-\eps^\ast|\cdot
\sum_{k\neq0}a(k)\sum_{i\in\Z}\big|f(x-e_i)-f(x)\big|.
\ec
\epro

We will show next that the position $l^\eps_t$ of the leftmost $1$ is
tight in the bias parameter $\eps$. First we need the following simple fact.

\bl[Stationary increments]
Let\label{L:statinc} $X$ and $Y$ be real random variables that are equal in
distribution, and assume that $\E\big[(Y-X)\vee 0]<\infty$. Then
$\E\big[|Y-X|\big]<\infty$ and $\E[Y-X]=0$.
\el
\bpro
It suffices to show that $\E\big[(X-Y)\vee 0\big]=\E\big[(Y-X)\vee
  0\big]$. For any real random variable $Z$ and constant $c>0$, write
$Z^c:=Z\wedge c$ and $Z_c:=Z\vee(-c)$. Then $\E[X^n_n-Y^n_n]=0$ and hence
$\E\big[(X^n_n-Y^n_n)\vee 0\big]=\E\big[(Y^n_n-X^n_n)\vee 0\big]$. By monotone
convergence
\be
\E\big[(X^n_n-Y^n_n)\vee 0\big]=\E\big[1_{\{-n<X\}}1_{\{Y<n\}}(X^n-Y_n)\vee 0\big]
\asto{n}\E\big[(X-Y)\vee 0\big],
\ee
and similarly $\E\big[(Y^n_n-X^n_n)\vee 0\big]$ converges to $\E\big[(Y-X)\vee
  0\big]$.
\epro

\bl[Expected displacement of the leftmost 1]
Let $l^\eps_t$ be the position of the leftmost 1 at time $t$ for the biased voter model $(X^\eps_t)_{t\geq0}$ with initial law $\nu_\eps$.
Then there\label{L:ltight} exists a constant $C<\infty$ such that uniformly in $\eps\in[0,1)$ and $t\geq 0$,
\be
\E\big[|l^\eps_t|\big]\leq Ct.
\ee
\el
\bpro
We first lower bound $l^\eps_t$. Since for any $\eps\in[0,1)$, the rate that
$(l^\eps_t)_{t\geq0}$ jumps to the left by $k$ ($k>0$) is given by
$\sum_{n\geq0}1_{\{X^\eps_t(l_t^\eps+n)=1\}}a(-n-k)$, we can couple it with a
uni-directional random walk $(S_t)_{t\geq0}$ started in $S_0=0$ with increment
rate $q(-k):=\sum_{n\geq0}a(-n-k)$ and $q(k)=0$ for all $k>0$, such
that $S_t\leq l^\eps_t$ for all $t\geq0$ almost surely. This gives the estimate
\be\label{leftinc}
\E\big[(-l^\eps_t)\vee 0\big]\leq\E[|S_t|]
=t\sum_{k\geq 1}k\sum_{n\geq0}a(-n-k)=t\sum_{k\geq 1}a(-k)\ha k(k+1).
\ee
The same argument applied to the rightmost zero gives
\be\label{rightinc}
\E\big[(r^\eps_t-r^\eps_0)\vee 0\big]\leq t\sum_{k\geq 1}a(k)\ha k(k+1).
\ee
Together, these bounds show that
\be
\E\big[\big\{(r^\eps_t-l^\eps_t)-r^\eps_0\big\}\vee 0\big]<\infty,
\ee
so we can apply Lemma~\ref{L:statinc} to the equally distributed random
variables $(r^\eps_t-l^\eps_t)$ and $r^\eps_0$ to conclude that
\be\label{stab}
\E\big[(r^\eps_t-l^\eps_t)-r^\eps_0\big]=0.
\ee
The bounds (\ref{leftinc}) and (\ref{rightinc}) show that $\E[l^\eps_t]$ is
well-defined (may be $+\infty$) and so is $\E[r^\eps_t-r^\eps_0]$ (may be $-\infty$). Therefore (\ref{stab}) implies
$\E[l^\eps_t]=\E\big[r^\eps_t-r^\eps_0\big]$. This gives
\be\ba{l}
\dis\E\big[l^\eps_t\vee 0\big]-\E\big[(-l^\eps_t)\vee 0\big]
=\E[l^\eps_t]=\E[r^\eps_t-r^\eps_0]\\[5pt]
\dis\quad
=\E\big[(r^\eps_t-r^\eps_0)\vee 0\big]-\E\big[(r^\eps_0-r^\eps_t)\vee 0\big]
\leq\E\big[(r^\eps_t-r^\eps_0)\vee 0\big],
\ec
and hence
\be
\E\big[|l^\eps_t|\big]
=\E\big[l^\eps_t\vee 0\big]+\E\big[(-l^\eps_t)\vee 0\big]
\leq\E\big[(r^\eps_t-r^\eps_0)\vee 0\big]+2\E\big[(-l^\eps_t)\vee 0\big]
\leq Ct,
\ee
with $C=\sum_{k\geq 1}\big[a(-k)+\ha a(k)\big]k(k+1)$.
\epro
\medskip

\bpro[ of Theorem~\ref{T:convpt}]
Let $\eps_n,\eps^\ast\in[0,1)$ and $\eps_n\to\eps^\ast$. We need to show that
$\nu_{\eps_n}\Rightarrow\nu_{\eps^\ast}$. Since $\{0,1\}^\Z$ is
compact under the product topology, $(\nu_{\eps_n})_{n\in\N}$ is tight. By
going to a subsequence if necessary, we may assume that $\nu_{\eps_n}\Rightarrow\nu^\ast$ for some
probability measure $\nu^\ast$ on $\{0,1\}^\Z$. 
If random variables $X^{\eps_n}$ and $X^*$ have law $\nu_{\eps_n}$ and $\nu^*$, respectively, then by Lemma~\ref{L:1bdbd}, $\sup_{n}\E[I_1(X^{\eps_n})]\leq C$ for some constant $C$.
Passing to the limit and applying Fatou's lemma gives $\E[I_1(X^*)]\leq C$,
which implies that $\nu^\ast$ is concentrated on $\{ x\in S^{01}_{\rm int} \cup S^{00}_{\rm int}:
x(0)=1\}$, (see~\eqref{eq:S01} for the definition of $S_{\rm int}^{00}$).
%{\color{blue} By Lemma~\ref{L:1bdbd} and Fatou's Lemma,
%$\nu^\ast$ is concentrated on $\{ x\in S^{01}_{\rm int} \cup S^{00}_{\rm int}:
%  x(0)=1\}$, see~\eqref{eq:S01}.} \Swa{[Like before, this does not seem to be
%    a typical application of Fatou. We should probably add some detail,
%    1.\ saying that if $X^\ast$ has law $\nu^\ast$, then $\E[I_1(X)]\leq C$ with
%    the same constant $C$ as in Lemma~\ref{L:1bdbd}, 2.\ pointing out this
%    implies $\nu^\ast$ is concentrated on $S^{01}_{\rm int}\cup S^{00}_{\rm
%      int}$, 3.\ and also explaining how 1.\ can be obtained from Fatou.]}
Let $X^{\eps_n}_\cdot$
denote the biased voter model with bias $\eps_n$ and initial law $\nu_{\eps_n}$. Fix $t>0$.
By Lemma~\ref{L:ltight}, $(l^{\eps_n}_t)_{n\geq 1}$ are tight, and hence by going to a subsequence if necessary, we
may assume that $(X^{\eps_n}_t,l^{\eps_n}_t)$ converges in law to some random variable $(X^\ast_t,l^\ast_t)$, 
which implies that 
\be\label{eq:nuconv}
\tet_{l^{\eps_n}_t}\big(X^{\eps_n}_t\big)\Asto{n}
\tet_{l^\ast_t}\big(X^\ast_t\big) \quad \mbox{and} \quad l^\ast_t=\min\{i:X^\ast_t(i)=1\}.
\ee
This can be seen by applying Skorohod's representation theorem
so that $(X^{\eps_n}_t,l^{\eps_n}_t)$ and $(X^*_t, l^*_t)$ are coupled in such a way that $(X^{\eps_n}_t,l^{\eps_n}_t)\to (X^*_t, l^*_t)$
almost surely. The left hand side of \eqref{eq:nuconv} has law $\nu_{\eps_n}$, and hence $\tet_{l^\ast_t}\big(X^\ast_t\big)$
has law $\nu^*$. On the other hand, Lemma~\ref{L:bvtcon} implies that $X^\ast_t$ is
distributed as the bias voter model with bias $\eps^\ast$, initial law $\nu^\ast$,
and evaluated at time $t$. Since $t>0$ is arbitrary, $\nu^*$ is an invariant law for the biased voter model with bias parameter $\eps^*$, seen from the 
leftmost 1. We note that the restriction of $\nu^*$ to $S_{\rm int}^{01}$ and $S_{\rm int}^{00}$ leads to two invariant measures for the biased voter model seen from the leftmost 1. However, it is not hard to see that the only such invariant measure on $S_{\rm int}^{00}$ is the zero measure. Therefore $\nu^*$ must concentrate on $S_{\rm int}^{01}$, and hence we must have $\nu^*=\nu_{\eps^*}$ by the uniqueness of the invariant law of the biased voter model on $\ov{S}_{\rm int}^{01}$.
%The only such invariant probability law on $\{ x\in S^{01}_{\rm
%  int} \cup S^{00}_{\rm int}: x(0)=1\}$ is $\nu_{\eps^*}$, because any such
%invariant law must assign full probability to $S^{01}_{\rm int}$ to avoid
%absorption in the state $x(\cdot)\equiv 0$. \Swa{[I would reformulate a
%    bit. We first argue that $\nu^\ast$ is concentrated in $S^{01}_{\rm
%      int}$. Then we use that the process on $S^{01}_{\rm int}$ has a unique
%    invariant law. The last step is rather hidden in the present formulation.]}
% Therefore we must
%have $\nu^*= \nu_{\eps^*}$.
\epro

\subsection{Equilibrium equation}\label{S:equi}

In this subsection, we will only consider the unbiased voter model. The
main purpose is to establish the equilibrium equation (\ref{equi}) below,
which shows that (\ref{equibdd}) holds with equality if $\eps=0$.  For
brevity, throughout this section we drop the superscript indicating the bias,
for example, the generator $G^0$ is abbreviated by $G$. Recall that %\sout{Let}
$X_\infty$ is the
%\sout{be a}
random variable taking values in $\hat S^{01}_{\rm int}$ (see
(\ref{hatS})) with law $\nu_0$, i.e., the law of the equilibrium voter interface
with the leftmost 1 shifted to the origin. 

\bp[Equilibrium equation]\label{P:equi}
Assume that the kernel $a(\cdot)$ has finite second moment $\sig^2=\sum_{k}a(k)k^2<\infty$.
Then, the following equilibrium equation for the voter model holds,
\be\label{equi2}
\E\big[Gh(X_\infty)\big]=0,
\ee
where $h=h^0$ is the number of inversions \eqref{numinv}.
Or equivalently, by the expression of $Gh$ in \eqref{Gh}, we have
\be\label{equi}
\E\Big[\sum_{k=1}^{\infty}a_{\rm s}(k)I_k(X_\infty)\Big]=\ha\sig^2,
\ee
where $a_{\rm s}(k)=\ha\big(a(k)+a(-k)\big)$ and $I_k(x)=|\{i:x(i)\neq x(i+k)\}|$.
\ep

We briefly explain our proof strategy.
If $h(X_t)-\int_{0}^{t}Gh(X_s)\di s$ was a martingale for any deterministic initial configuration $X_0$,
and if $\E[h(X_\infty)]$ and $\E[|Gh(X_\infty)|]$ were finite, then the equilibrium equation (\ref{equi2}) would follow.
However, two difficulties arise in this approach.
As shown in Remark~\ref{R:nonexp}, we can only prove that $h(X_t)-\int_{0}^{t}Gh(X_s)\di s$ is a martingale when $a(\cdot)$ has finite third moment, as otherwise there is no control on the expected number of inversions $\E[h(X_t)]$.
Worse still, since $h$ is bounded from below by the length $L$ of the interface and $\E[L(X_\infty)]=\infty$ by \cite[Theorem~6]{CD95} or \cite[Theorem~1.4]{BMSV06}, we have $\E[h(X_\infty)]=\infty$.
To bypass these difficulties, we will show that the equilibrium equation (\ref{equi2}) holds for $h_M$ in place of $h$, where
\be
h_M(x)=\big|\{(j,i):0<i-j\leq M,~x(j,i)=10\}\big|,
\ee
which serves as a truncation approximation of $h(x)$.
We will then let $M\to\infty$ to deduce (\ref{equi2}).
To deduce this last convergence, we will use the fact that the expected number
of 1-boundaries $\E[I_1(X_\infty)]$ is finite, by Lemma~\ref{L:1bdbd}.

%If $a$ has finite third moment, then as proved in Remark~\ref{R:nonexp}, we have that
%$h(X_t)-\int_{0}^{t}Gh(X_s)\di s$ is a martingale whenever
%$X_0$ is a deterministic configuration.
%\weg{Indeed, we can moreover show that this is a martingale even if
%  $(X_t)_{t\geq0}$ starts from the equilibrium $X_\infty$, and hence}
%We will moreover show that $\E[h(X_\infty)]$ and $\E[|Gh(X_\infty)|]$ are
% finite. Combining this with the martingale property, the equilibrium equation (\ref{equi2}) follows.
%To prove (\ref{equi2}) assuming only a finite second moment on $a$, the basic idea is to approximate the number of inversions $h(x)$ by a truncated version $h_M(x)=\big|\{(i,j):0<j-i\leq M,~x(i,j)=10\}\big|$ as defined in (\ref{htrun}), and then pass to the limit.

Our first step is to do a generator calculation for $h_M$.
Recall formula (\ref{G0}) for $Gh$, where $h$ is the number of inversions.
%\be
%Gh(x)=\sum_{k=1}^{\infty}a_{\rm s}(k)\big(k^2-I_k(x)\big).
%\ee
The next lemma identifies $Gh_M$.

\bl\label{L:htrun}
For any $x\in S^{01}_{\rm int}$ and $M\in\N$, we have
\be\label{Ghhat0}
Gh_M(x)=\sum_{k=1}^{\infty}a_{\rm s}(k)\big(k^2-I_k(x)\big)
+\sum_{n=1}^{\infty}A(n)I_{M+n}^{10}(x)-\sum_{n=1}^{\infty}A(n)I^{10}_{M-1-n}(x),
\ee
where
$A(n):=\dis\sum_{k=n}^{\infty}\big(a(k)+a(-k)\big)=2\sum_{k=n}^{\infty}a_{\rm s}(k)$.
\el
\bpro
We use the generator decomposition $G=\sum_{k\neq0}a(k)G_k$ in (\ref{Gksum}), and separately calculate $G_kh_M$ for $k>0$ and $k<0$.

For $k>0$, to calculate $G_kh_M(x)$, we consider all triples $(i,n,n-k)$ with
$|i-n|\leq M$, where an inversion in $x(i,n)$ is either created or destroyed because $x(n)$ changes its value to that of $x(n-k)$. Therefore,
\bc\label{Ghhat}
\dis G_kh_M(x)&=&\dis\sum_{i}1_{\{x(i)=1\}}\Big\{-\sum_{n=i+1}^{i+M}1_{\{x(n-k,n)=10\}}
+\sum_{n=i+1}^{i+M}1_{\{x(n-k,n)=01\}}\Big\}\\ [5pt]
&&\dis +\sum_{i}1_{\{x(i)=0\}}\Big\{\sum_{n=i-M}^{i-1}1_{\{x(n-k,n)=10\}}
-\sum_{n=i-M}^{i-1}1_{\{x(n-k,n)=01\}}\Big\}\\ [5pt]
&=&\dis\sum_{i}1_{\{x(i)=1\}}\sum_{n=i+1}^{i+M}\big(1_{\{x(n-k)=0\}}-1_{\{x(n)=0\}}\big)\\ [5pt]
&&+\dis\sum_{i}1_{\{x(i)=0\}}\sum_{n=i-M}^{i-1}\big(1_{\{x(n-k)=1\}}-1_{\{x(n)=1\}}\big),
\ec
where in the last equality we used (\ref{reexp}).
To further simplify this,
we apply the summation identity (\ref{rearr}) with $a=i+1-k,b_1=i+1,b_2=i+M-k$ and $c=i+M$, which gives
\be
\sum_{n=i+1}^{i+M}1_{\{x(n)=0\}}-\sum_{n=i+1}^{i+M}1_{\{x(n-k)=0\}}
=\sum_{n=-k+1}^{0}1_{\{x(i+M+n)=0\}}-\sum_{n=-k+1}^{0}1_{\{x(i+n)=0\}},
\ee
or equivalently,
\be
\sum_{n=i+1}^{i+M}\big(1_{\{x(n-k)=0\}}-1_{\{x(n)=0\}}\big)
=\sum_{n=0}^{k-1}\big(1_{\{x(i-n)=0\}}-1_{\{x(i+M-n)=0\}}\big).
\ee
Similarly, we can also get
\be
\sum_{n=i-M}^{i-1}\big(1_{\{x(n-k)=1\}}-1_{\{x(n)=1\}}\big)
=\sum_{n=1}^k\big(1_{\{x(i-M-n)=1\}}-1_{\{x(i-n)=1\}}\big).
\ee
Substituting the above identities into the right-hand side of (\ref{Ghhat}), and then using
the notation $I^{01}_k,I^{10}_k$ and $I_k$ as in (\ref{Idef}), we can rewrite
(\ref{Ghhat}) as
\bc\label{Ghhat2}
\dis G_kh_M(x)&=&\dis\sum_i\sum_{n=0}^{k-1}
\big(1_{\{x(i-n,i)=01\}}-1_{\{x(i,i+M-n)=10\}}\big)\\[5pt]
&&\dis+\sum_i\sum_{n=1}^k\big(1_{\{x(i-M-n,i)=10\}}-1_{\{x(i-n,i)=10\}}\big)\\[5pt]
&=&\dis\sum_{n=1}^{k-1}I^{01}_n(x)-\sum_{n=0}^{k-1}I^{10}_{M-n}(x)
-\sum_{n=1}^{k}I^{10}_{n}(x)+\sum_{n=1}^{k}I^{10}_{M+n}(x)\\ [5pt]
&=&\dis\ha\big(k^2-I_k(x)\big)
+\sum_{n=1}^{k}\big(I^{10}_{M+n}(x)-I^{10}_{M-1-n}(x)\big).
\ec
where in the last equality we applied (\ref{Irel}) to
$\sum_{n=1}^{k-1}\big(I^{01}_n(x)-I^{10}_n(x)\big)-I^{10}_k(x)$.

Using symmetry, we can %\sout{now}
also easily obtain a formula for $G_kh_M$ when
  $k<0$. For any $x\in S^{01}_{\rm int}$, define $x'\in S^{01}_{\rm int}$ by
  $x'(i):=1-x(-i)$ $(i\in\Z)$. Then, for any function $f:S^{01}_{\rm
    int}\to\R$, one has $G_kf(x)=G_{-k}f'(x')$, where $f'(x):=f(x')$ $(x\in
  S^{01}_{\rm int})$. We observe that $I_k(x)=I_k(x')$ and hence also
  $h_M(x)=\sum_{k=1}^MI_k(x)$ is symmetric in the sense that $h_M(x)=h_M(x')$
  $(x\in S^{01}_{\rm int})$. Combining these observations with (\ref{Ghhat2}),
  we obtain that for any $k\neq 0$
\be\label{Ghhat3}
G_kh_M(x)=\ha\big(k^2-I_{|k|}(x)\big)
+\sum_{n=1}^{|k|}\big(I^{10}_{M+n}(x)-I^{10}_{M-1-n}(x)\big).
\ee
Inserting this into $G=\sum_{k\neq 0}a(k)G_k$,
we have
\be\label{Ghhat4}
Gh_M(x)
=\sum_{k=1}^{\infty}a_{\rm s}(k)\big(k^2-I_k(x)\big)
+2\sum_{k=1}^{\infty}a_{\rm s}(k)\sum_{n=1}^{k}
\big(I^{10}_{M+n}(x)-I^{10}_{M-1-n}(x)\big).
\ee
%\bc\label{Ghhat4}
%\dis Gh_M(x)%&=&\dis\sum_{k\neq0}^{\infty}a(k)G_kh_M(x)\\ [5pt]
%&=&\dis\sum_{k=1}^{\infty}a_{\rm s}(k)\big(k^2-I_k(x)\big)
%+2\sum_{k=1}^{\infty}a_{\rm s}(k)\sum_{n=1}^{k}
%\big(I^{10}_{M+n}(x)-I^{10}_{M-1-n}(x)\big)\\ [5pt]
%&=&\dis\sum_{k=1}^{\infty}a_{\rm s}(k)\big(k^2-I_k(x)\big)
%+\sum_{n=1}^{\infty}A(n)\big(I^{10}_{M+n}(x)-I^{10}_{M-1-n}(x),
%\ec
Interchanging the summation order, we obtain (\ref{Ghhat0}).
\epro

We are now ready to prove Proposition~\ref{P:equi}.\med

\bpro[~of Proposition \ref{P:equi}]Let $(X_t)_{t\geq0}$ be the voter model starting from some deterministic configuration $x\in S^{01}_{\rm int}$.
Under the second moment assumption, by %\sout{(\ref{hMmart})}
Lemma~\ref{L:hMmart}, for any $t>0$,
\be\label{equi5}
\E\big[h_M(X_t)\big]-h_M(x)=
\int_{0}^{t}\E\big[G h_M(X_s)\big]\di s.
\ee
Assume for the moment that both $h_M$ and $Gh_M$ are absolutely integrable with respect to the law of $X_\infty$. Then we can integrate both sides of \eqref{equi5} with respect to the invariant law to get
\be\label{equi3}
\E\big[Gh_M(X_\infty)\big]=0.
\ee
By letting $M\to\infty$, we will see in the following that \eqref{equi3} implies \eqref{equi2}, the equilibrium equation for the voter model.
Recalling the expression of $Gh_M(x)$ in (\ref{Ghhat0}) and $Gh(x)$ in (\ref{G0}), we obtain from (\ref{equi3}) that
\be\label{equi4}
\ba{l}
\dis\E\big[Gh(X_\infty)\big]=
\dis\E\Big[\sum_{k=1}^{\infty}a_{\rm s}(k)\big(k^2-I_k(X_\infty)\big)\Big]\\ [5pt]
\qquad\dis=\E\Big[\sum_{n=1}^{\infty}A(n)I^{10}_{M-1-n}(X_\infty) -\sum_{n=1}^{\infty}A(n)I^{10}_{M+n}(X_\infty)\Big],
\ec
where $A(n)=2\sum_{k=n}^{\infty}a_{\rm s}(k)$.
For $n\geq M$, by %\sout{(\ref{Idef})}
(\ref{Idef2}) and (\ref{Irel}), we have
\be
I^{10}_{-(n+1-M)}(x)=I^{01}_{n+1-M}(x)=I^{10}_{n+1-M}(x)+(n+1-M).
\ee
Therefore, by Lemma \ref{L:Ibdd}, we can bound the difference in the expectation in \eqref{equi4} by
\be\label{equi6}
\ba{l}
\dis\big|\sum_{n=1}^{\infty}A(n)I^{10}_{M-1-n}(X_\infty) -\sum_{n=1}^{\infty}A(n)I^{10}_{M+n}(X_\infty)\big|\\ [5pt]
\quad\leq\dis\sum_{n=1}^{M-1}A(n)\big|I^{10}_{M-1-n}(X_\infty)-I^{10}_{M+n}(X_\infty)\big|\\ [5pt]
\qquad\dis+\sum_{n=M}^{\infty}A(n)\big|I^{10}_{n+1-M}(X_\infty)+(n+1-M)-I^{10}_{M+n}(X_\infty)\big|\\ [5pt]
\quad\leq\dis\sum_{n=1}^{M-1}A(n)(2n+1)I^{10}_1(X_\infty)+
\sum_{n=M}^{\infty}A(n)\big\{n+1-M+(2M-1)I^{10}_1(X_\infty)\big\}\\ [5pt]
%\quad\leq\dis\sum_{n=1}^\infty A(n)(2n+1)I^{10}_1(X_\infty)+
%\sum_{n=0}^{\infty}A(n+M)(n+1)\\ [5pt]
\quad\leq\dis\sum_{n=1}^{\infty}3nA(n) \big(1+I^{10}_1(X_{\infty})\big).
\ec
Due to the second moment assumption, $\sum_{n=1}^{\infty}nA(n)$ is finite,
so applying Lemma~\ref{L:1bdbd} we see that the right-hand side of
(\ref{equi6}) is bounded in expectation. As a result, the term in the
expectation on the right-hand side of (\ref{equi4}) is bounded in
  absolute value by an integrable random variable, uniformly in $M$.
If this term moreover converges to zero pointwise as $M$ tends to $\infty$, then
applying the dominated convergence theorem to (\ref{equi4}), we will obtain the equilibrium equation
\be
\E\big[Gh(X_\infty)\big]=0.
\ee
To see the pointwise convergence, note that for every $x\in S^{01}_{\rm int}$, there exists some $M_x$ such that $I_k(x)=0$ for all $|k|> M_x$, and thus when $M>M_x+1$,
\be
\sum_{n=1}^{\infty}A(n)I^{10}_{M-1-n}(x)-\sum_{n=1}^{\infty}A(n)I^{10}_{M+n}(x)=\sum_{k=-M_x}^{M_x}A(M-1-k)I_k(x),
\ee
where the right-hand side decreases to $0$ since $A(M-1-k)\down0$ as $M$ tends to $\infty$.

To complete the proof, it remains to show, for fixed $M$, the absolute integrability of $h_M$ and $Gh_M$ with respect to the invariant law.
For the nonnegative function $h_M$, by Lemmas~\ref{L:Ibdd} and \ref{L:1bdbd},
\be
\E\big[h_M(X_\infty)\big]=\E\Big[\sum_{k=1}^{M}I^{10}_k(X_\infty)\Big]\leq\E\Big[\sum_{k=1}^{M}k I^{10}_1(X_\infty)\Big]<\infty.
\ee
By using the expression of $Gh_M$ in \eqref{Ghhat0}, Lemma~\ref{L:Ibdd},
the fact that $\sum_{n=1}^\infty A(n)n<\infty$ since $\sum_ka(k)k^2<\infty$,
and Lemma~\ref{L:1bdbd}, it is also not hard to see that
$\E\big[\big|Gh_M(X_\infty)\big|\big]<\infty$.
\epro

\subsection{Proof of Theorem~\ref{T:nucont}}\label{S:ConS}
%{\jan [Below has been rewritten. In the original formulation it was claimed
%    that we prove tightness of the family $(\nu_\eps)_{\eps\in(0,\ga)}$ for
%    some $\ga\in(0,1)$. This is incorrect. The error occurred at the place
%    where we said:} ``We begin by assuming that
%  $\big(L(X^\eps_0)\big)_{\eps\in(0,\ga)}$ is not tight.  Then there exists
%  $\delta>0$ such that for all $N>0$, there exists $\eps_N$
%  {\color{red}($\eps_N\down0$ as $N\to\infty$)}''. {\jan The claim in red is
%    unsubstantiated. It is possible that
%    $\big(L(X^\eps_0)\big)_{\eps\in(0,\ga)}$ is not tight yet we cannot choose
%    the $\eps_N$ in such a way that $\eps_N\down0$.]}\med

For each $\eps\geq 0$, let $X^\eps_\infty$ denote a random variable taking
values in $\hat S^{01}_{\rm int}$ (see (\ref{hatS})) with law
$\nu_\eps$. It suffices to show that as $\eps\down 0$, the measures $\nu_\eps$ converge weakly to $\nu_0$ with respect
to the discrete topology on $\hat S^{01}_{\rm int}$. By
Theorem~\ref{T:convpt}, the measures $\nu_\eps$ converge weakly to $\nu_0$
with respect to the product topology on $\{0,1\}^\Z$. To improve this to
convergence with respect to the discrete topology, it suffices to show that for any sequence
$\eps_n\down 0$, the laws of the random variables $\big(r^{\eps_n}_\infty\big)_{n\geq 1}$
are tight, where as in (\ref{left}), we let
$r^\eps_\infty:=\max\{i:X^\eps_\infty(i)=0\}$ denote the position of the
rightmost zero of $X^\eps_\infty$.

Assume that for some $\eps_n\down 0$, the laws of $\big(r^{\eps_n}_\infty\big)_{n\geq 1}$ are not tight. Then going to a subsequence if necessary, we
can find $\de>0$ and $(m_N)_{N\geq 1}$ such that
\be\label{deN}
\P\big[r^{\eps_n}_\infty>N\big]>\de\qquad \mbox{for all  } n\geq m_N.
\ee
For $x\in S^{01}_{\rm int}$ and $N\in\Z$, let
\be
I^N_k(x):=\big|\{i\leq N:x(i)\neq x(i+k)\}\big|.
\ee
Since $X^\eps_\infty(r^\eps_\infty)=0\neq1=X^\eps_\infty(r^\eps_\infty+k)$ for
all $k\geq 1$, by Lemma~\ref{L:equibdd} and (\ref{deN}), we see that
\be
\frac{1}{2} \sig^2 \geq \E\big[\sum_{k=1}^{\infty}a_{\rm s}(k)I_k(X^{\eps_n}_\infty)\big] \geq
\E\big[\sum_{k=1}^{\infty}a_{\rm s}(k)I^N_k(X^{\eps_n}_\infty)\big]
+A\de\qquad \mbox{for all  } n\geq m_N,
\ee
where $A=\sum_{k=1}^\infty a_{\rm s}(k)>0$. Letting $n\to\infty$ and using the weak convergence of $X^{\eps_n}_\infty$ to $X^0_\infty$ with respect to the product topology on $\{0,1\}^\Z$ and Fatou's Lemma, we find that
\be
\E\big[\sum_{k=1}^{\infty}a_{\rm s}(k)I^N_k(X^0_\infty)\big]
\leq\ha\sig^2-A\de.
\ee
As we send $N\to\infty$, the left hand side converges as $I^N_k\up I_k$, which leads to a contradiction with (\ref{equi}) in Proposition~\ref{P:equi}.\QED

\med\med

\noi
{\bf Acknowledgement} J.M.~Swart is sponsored by grant 16-15238S of the Czech
Science Foundation (GA CR). R.~Sun and J.~Yu are partially supported by NUS
grant R-146-000-220-112. We would like to thank the Institute for
Mathematical Sciences at NUS for hospitality during the program
\emph{Genealogies of Interacting Particle Systems}, where part of this work was
done.
We also thank an anonymous referee for a careful reading of the paper and helpful suggestions.

%Note: dedication changed, no longer to 15-08819S.

%\Swa{[Changed the order of the first two entries in the references.]}

\end{document}